\documentclass[12pt]{article}
\usepackage{amsmath,amsthm,amsfonts,latexsym,amssymb,fancyhea,supertab}
\usepackage{a4}

  \def\Q{\mathbb{Q}}
\def\Z{\mathbb{Z}} \def\R{\mathbb{R}}

\def\cE{{\cal E}}  \def\cF{{\cal F}}

\def\Nm{\mathrm{Norm}}

 \def\be{\beta} \def\ep{\epsilon} 
\def\th{\theta} \def\ph{\phi} \def\ga{\gamma}

\def\vn{\vspace{1ex}\noindent}

\def\Proof{\noindent{\bf Proof.  }}

\newcommand{\proofend}{\hspace*{1mm} \hfill{$\Box$}}

\begin{document}

\newtheorem*{u8_result}{Theorem}
\newtheorem{problem}{Problem}
\title{Lucas sequences whose 8th term is a square}
\author{A.~Bremner\thanks{Department of Mathematics, Arizona State University, Tempe AZ, USA, e-mail: bremner@asu.edu, http://\~{}andrew/bremner.html} \and 
N.~Tzanakis\thanks{Department of Mathematics, University of Crete,
Iraklion, Greece, e-mail: tzanakis@math.uoc.gr , http://www.math.uoc.gr/\~{}tzanakis}
}
\maketitle
\section{Abstract}
Let $P$ and $Q$ be non-zero integers. The Lucas sequence $\{U_n(P,Q)\}$
is defined by
\[ U_0=0, \quad U_1=1, \quad U_n= P U_{n-1}-Q U_{n-2} \quad (n \geq 2). \]
For each positive integer $n \leq 7$ we describe all Lucas sequences
with $(P,Q)=1$ having the property that $U_n(P,Q)$ is a perfect square.
The arguments are
elementary. We also find all Lucas sequences such that $U_8(P,Q)$
is a perfect square. This reduces to a number of problems of similar type,
namely, finding all points on an elliptic curve defined over a quartic
number field subject to a ``{\bf Q}-rationality'' condition on
the $X$-coordinate. This is achieved by $p$-adic computations (for a
suitable prime $p$) using the formal group of the elliptic curve.
\section{Introduction} \label{introduction}
Let $P$ and $Q$ be non-zero integers. The Lucas sequence $\{U_n(P,Q)\}$ 
is defined by
\begin{equation}
\label{Lucas}
U_0=0, \quad U_1=1, \quad U_n= P U_{n-1}-Q U_{n-2} \quad (n \geq 2).
\end{equation}
The sequence $\{U_n(1,-1)\}$ is the familiar Fibonacci sequence, and it was proved by 
Cohn~\cite{Co1} in 1964 that the only perfect square
greater than $1$ in this sequence is $U_{12}=144$. The question arises, 
for which parameters $P$, $Q$, can $U_n(P,Q)$ be a perfect square? 
In what follows, we shall assume that we are not dealing with the degenerate 
sequences corresponding to $(P,Q)=(\pm 1,1)$, where $U_n$ is periodic
with period $3$, and we also assume $(P,Q) \neq (-2,1)$ (in which case 
$U_n=\Box$ precisely when $n$ is an odd square) and  $(P,Q) \neq (2,1)$ 
(when $U_n=\Box$ precisely when $n$ is square).
Ribenboim and McDaniel~\cite{RM1} with only elementary methods show that when $P$ and $Q$ are 
{\it odd}, and $P^2-4Q>0$,
then $U_n$ can be square only for $n=0,1,2,3,6$ or $12$; and that there are at most two indices 
greater than 1 for which $U_n$ can be square.
They characterize fully the instances when $U_n=\Box$, for $n=2,3,6$.
Bremner \& Tzanakis~\cite{BT} extend these results by determining all Lucas 
sequences $\{U_n(P,Q)\}$ with $U_{12}=\Box$, subject only to the restriction 
that gcd$(P,Q)=1$ (it turns
out that the Fibonacci sequence provides the only example). Under the same hypothesis, all 
Lucas sequences with $\{U_n(P,Q)\}$ with $U_9=\Box$ are determined.
There seems little mention in the literature of when
under general hypotheses $U_n(P,Q)$ can be a perfect square. Note that 
for $n \geq 1$, $U_n(k P, k^2 Q) = k^{n-1} U_n(P,Q)$, and so for fixed $P$, $Q$, 
and {\it even} $n$, appropriate choice of $k$ gives a sequence
with $U_n(k P, k^2 Q)$ a perfect square. The restriction to $(P,Q)=1$ is therefore
a sensible one, and we shall assume this from now on.  A small computer search reveals 
sequences with $U_n(P,Q)$ a perfect square, only for $n=0,\ldots,8$, and $n=12$. 
Bremner \& Tzanakis~\cite{BT} have addressed the case $n=12$. Section \ref{up to U7}
of this paper addresses the case of $U_n(P,Q)=\Box, \; n\leq 7$,
which can be treated entirely elementarily. The remainder of the paper (section \ref{case U8})
addresses the case $U_8(P,Q)=\Box$.
This reduces to a number of problems of similar type, namely, finding all points on an 
elliptic curve defined over a number field $K$ subject to a ``$\Q$-rationality" condition
on the $X$-coordinate. The elliptic curves we consider have $K$-rank at most 2, with
degree $[K:\Q]=4$, and so
this problem is of ``Chabauty" type in the language of Nils Bruin. Bruin has powerful
techniques for addressing this type of problem, and~\cite{Br1},~\cite{Br2},~\cite{Br3},~\cite{Br4} provide details and examples. We
persevere in writing the current paper to describe in very concrete form the underlying mathematics, 
based on the work of Flynn and Wetherell~\cite{FW}, together with a theorem that is 
essentially due to Th.~Skolem from the 1930s to deal with the example of our rank 2 
elliptic curve. The latest release of Magma now contains Bruin's routines for much of 
the calculations of this paper, but we feel it is still worthwhile to give some (minimal) 
details of the computations, in order to expose the underlying theory and make it 
accessible to the reader, as well as for those without access to Magma. 
\section{Solution of $U_n(P,Q)=\Box, \; n\leq 7$} \label{up to U7}
Certainly $U_2(P,Q)=\Box$ if and only if $P=a^2$, and $U_3(P,Q)=\Box$ if and
only if $P^2-Q=a^2$.\\
Now $U_4(P,Q)=\Box$ if and only if $P(P^2-2Q)=\Box$, so if and only if either 
$P= \delta a^2, Q=\frac{1}{2}(a^4-\delta b^2)$, or $P=2\delta a^2, Q=2a^4-\delta b^2$,
with $\delta=\pm 1$ (where, in the first instance, $a b$ is odd and in the
second instance $b$ is odd). \\
The demand that $U_5(P,Q)$ be square is that $P^4-3 P^2 Q+Q^2 = \Box$, equivalently,
that $1-3 x+x^2=\Box$, where $x=Q/P^2$. Parametrizing the quadric,
$Q/P^2 = (5 \lambda^2+6 \lambda \mu +\mu^2)/(4 \lambda \mu)$, where, without loss
of generality, $(\lambda, \mu)=1$, $\lambda>0$, and $\mu \not \equiv 0 \pmod 5$. Necessarily
$(\lambda, \mu)=(a^2, \pm b^2)$, giving $(P,Q)=(2 a b, 5 a^4 + 6 a^2 b^2 +b^4)$ or
$(2 a b, -5 a^4+6 a^2 b^2 - b^4)$  if $a$ and $b$ are of opposite parity, and 
$(P,Q)=(a b, \frac{1}{4}(5 a^4 + 6 a^2 b^2 +b^4))$ or 
$(a b, \frac{1}{4}(-5 a^4 + 6 a^2 b^2 - b^4))$, if $a$ and $b$ are both odd.\\
The demand that $U_6(P,Q)$ be square is that $P(P^2-Q)(P^2-3Q)=\Box$, which leads to
one of seven cases: $P=a^2$, $P^2-Q=b^2$, with $-2 a^4+3 b^2=\Box$; $P=a^2$, $P^2-Q=-2 b^2$,
with $a^4+3 b^2=\Box$; $P=-a^2$, $P^2-Q=2 b^2$, with $a^4-3 b^2=\Box$; and $P=3 a^2$,
$P^2-Q= \delta b^2$, ($\delta = \pm 1, \pm 2$), with $-\frac{6}{\delta} a^4+b^2=\Box$.
So finitely many parametrizations result (which can easily be obtained, if we wish
to do so).\\
The demand that $U_7(P,Q)$ be square is that $P^6-5P^4 Q+6 P^2 Q^2-Q^3 =\Box$,
equivalently, that $1+5 x+6 x^2 +x^3 = y^2$, where $x=-Q/P^2$. This latter elliptic curve
has rank 1, with generator $P_0=(-1,1)$, and trivial torsion. Accordingly,
sequences with $U_7(P,Q)=\Box$ are parametrized by the multiples of $P_0$ on the  above
elliptic curve, corresponding to $(\pm P,Q) = (1,1)$, $(1,5)$, $(2,-1)$, $(5,21)$, $(1,-104)$, 
$(21,545)$, $(52,415)$,...\\
\section{Solution of $U_8(P,Q)=\Box$} \label{case U8}
The remainder of the paper will be devoted to the proof of the following result: 
\begin{u8_result} \label{u8_result}
The only non-degenerate sequences where $(P,Q)=1$ and $U_8(P,Q)=\Box$ are given by 
$U_8(1,-4)=21^2$ and $U_8(4,-17)=620^2$.
\end{u8_result}
\subsection{The auxiliary equations} \label{auxiliary eqs}
The demand that $U_8(P,Q)$ be square is that $P(P^2-2Q)(P^4-4P^2Q+2Q^2)=\Box$. 
\subsubsection{$P$ odd} \label{odd P}
It follows that $(P,P^2-2Q,P^4-4P^2Q+2Q^2) = (a^2,b^2,c^2)$, $(a^2,-b^2,-c^2)$,
$(-a^2,b^2,-c^2)$, or $(-a^2,-b^2,c^2)$, where $a$, $b$, $c$ are positive integers 
with $a b$ odd. The latter two possibilities are impossible modulo $4$,
and the first two possibilities lead respectively to:
\begin{eqnarray}
-a^8+2 a^4 b^2 +b^4 & = & 2 c^2 \label{eq1} \\
-a^8-2 a^4 b^2 +b^4 & = & - 2 c^2 \label{eq2}
\end{eqnarray}
Equation (\ref{eq1}) is related to the elliptic curves ${\mathcal E}_1$ and ${\mathcal E}_2$ 
(see (\ref{ell1}) and (\ref{ell2}), respectively) and
equation (\ref{eq2}) is related to the elliptic curves ${\mathcal E}_3$ and ${\mathcal E}_4$ 
(see (\ref{ell3}) and (\ref{ell4}), respectively).
According to Proposition \ref{Results_r1} the only positive solutions to the above equations 
are $(a,b)=(1,1),(1,3)$ and $(1,1)$ respectively, leading
to $(P,Q)=(1,0), (1,-4)$ and $(1,1)$, from which we reject the first one.
The last gives a degenerate sequence.
\subsubsection{$P$ even} \label{even P}
Now $Q$ is odd, and $2$ exactly divides both $P^4-4P^2Q+2Q^2$ and $P^2-2Q$,
forcing $P \equiv 0 \pmod 4$. Put $P=4p$, so that $U_8=\Box$ if and only if
$p(8p^2-Q)(128p^4-32p^2Q+Q^2)=\Box$, with $(p,Q)=1$. It follows that
$(p,8p^2-Q,128p^4-32p^2Q+Q^2)=(a^2,b^2,c^2)$, $(a^2,-b^2,-c^2)$, $(-a^2,b^2,-c^2)$,
or $(-a^2,-b^2,c^2)$, where $a$, $b$, $c$ are positive integers, $(a,b)=1$ 
and $b c$ is odd. The middle two possibilities are impossible modulo $4$, 
and the remaining two possibilities lead respectively to:
\begin{eqnarray}
-64a^8+16a^4b^2+b^4 & = & c^2 \label{eq3} \\
-64a^8-16a^4b^2+b^4 & = & c^2 \label{eq4}
\end{eqnarray}
Equation (\ref{eq3}) is related to the elliptic curves 
${\mathcal E}_5,{\mathcal E}_6,{\mathcal E}_7$ and ${\mathcal E}_8$ 
(see (\ref{ell5}), (\ref{ell6}), (\ref{ell7}) and (\ref{ell8}), respectively).
According to Proposition \ref{Results_r1} the only positive solution which leads 
to a desired pair $(P,Q)$ is $(a,b)=(1,5)$, leading to $(P,Q)=(4,-17)$. \\
Equation (\ref{eq4}) is related to the elliptic curves 
${\mathcal E}_9,{\mathcal E}_{10},{\mathcal E}_{11}$ and ${\mathcal E}_{12}$ 
(see (\ref{ell9}), (\ref{ell10}), (\ref{ell11}) and (\ref{ell12}), respectively).
According to Proposition \ref{Results_r1}, which deals with ${\mathcal E}_i$ with $i=9,11,12$ 
and Proposition \ref{Results_r2}, which deals with 
${\mathcal E}_{10}$, there are no positive solutions $(a,b)$. 

\subsection{The elliptic curves} \label{elliptic curves}
In this section we reduce the solution of equations (\ref{eq1})-(\ref{eq4}) to the
solution of a number of problems all of which fit the following general shape: 
\begin{problem} \label{problem}
Let 
\begin{equation} \label{gen_xy_elliptic}
\cE\,:\, Y^2+a_1XY+a_3Y=X^3+a_2X^2+a_4X +a_6
\end{equation} 
be an elliptic curve defined over $\Q(\alpha)$, where $\alpha$ is a root of a 
polynomial $f(X)\in\Z[X]$, irreducible over $\Q$, of degree $d\geq 2$, 
and let $\beta,\gamma \in\Q(\alpha)$ be algebraic integers. Find all points 
$(X,Y)\in\cE(\Q(\alpha))$ for which $\beta X+\gamma$ is a {\em rational number}.
\end{problem}
We shall see that equations (\ref{eq1})-(\ref{eq4}) lead to elliptic curves ${\mathcal E}_i$, $i=1,\ldots,12$, 
and so 12 instances of Problem \ref{problem}; in each case we 
specify the corresponding ``condition on $X$-coordinate" $\be X+\ga \in \Q$.
In all but one case the elliptic curves have rank 1 and in the 
exceptional case the rank is 2. \\
We need details of two number fields.
First, let $\theta$ be a root of $f_1(x)=x^4+2 x^2-1$, with $K_1=\Q(\theta)$.
The class number of $K_1$ is 1, the maximal order $\mathcal{O}_1$ of $K_1$
is $\Z[\theta]$, and fundamental units of $\mathcal{O}_1$ are
$\eta_1=\theta$, $\eta_2=2-3\theta+\theta^2-\theta^3$.
The factorization of 2 is $2=\eta_1^{-4} \eta_2^2 (1+\theta)^4$.\\
Second, let $\phi$ be a root of $f_2(x)=x^4+4 x^2-4$, with $K_2=\Q(\phi)$.
The class number of $K_2$ is 1, the maximal order $\mathcal{O}_2$ is
$\Z[1,\phi,\frac{1}{2}\phi^2,\frac{1}{2}\phi+\frac{1}{4}\phi^3]$, and
fundamental units are $\epsilon_1=\frac{1}{2}\phi+\frac{1}{4}\phi^3$,
$\epsilon_2=2+2\phi+\frac{1}{2}\phi^2+\frac{1}{2}\phi^3$. The factorization
of $2$ is $2 = \epsilon_2^{-2} \pi^4$, where $\pi= 1+\frac{3}{2}\phi+\frac{1}{4}
\phi^3$.
\subsubsection{Equation (\ref{eq1}) and curves ${\mathcal E}_1,{\mathcal E}_2$} \label{curve1}
The factorization of (\ref{eq1}) over $K_1$ is
\[ (b-\theta a^2)(b+\theta a^2)(b^2+(2+\theta^2) a^4) = 2 \Box, \]
and it is easy to see that the gcd of any two (ideal) terms on the left hand side is equal to 
$(1+\theta)$, with the last term exactly divisible by $(1+\theta)^2$. Hence,
\[(b+\theta a^2) (b^2+(2+\theta^2)a^4) = \pm\eta_1^i\eta_2^j(1+\theta)\Box, \]
where $i,j = 0,1$.  Specializing $\theta$ at the real root $0.643594...$ of $f_1(x)$, and
using $b>0$, then necessarily the sign on the right hand side must be positive.
By putting $b/a^2 = \delta^{-1} x/(1+\theta)$, where $\delta = \eta_1^i \eta_2^j$, our problem 
reduces to finding all $K_1$-points $(x,y)$ on the curves
\[ (x+\theta (1+\theta) \delta)(x^2+(2+\theta^2)(1+\theta)^2 \delta^2) = y^2, \]
subject to $\delta^{-1}x/(1+\theta) \in \Q$, with $\delta=1$, $\eta_1$, $\eta_2$, or 
$\eta_1 \eta_2$.  Putting
\[ (x,y) = (2 X-(\theta+\theta^2) \delta, 2(1+\theta^2) Y) \]
gives
\begin{equation}
\label{E_1}
Y^2 = X( X^2 -(\theta+\theta^2) \delta X +(1+\theta+\theta^3) \delta^2);
\end{equation}
and the condition on the $X$-coordinate becomes:
\[ -\theta+\frac{(3-3\theta+\theta^2-\theta^3) X}{\delta} \in \Q. \]
There are several computer packages now available for computing with elliptic
curves $E$ over number fields $K$. We mention Algae~\cite{BrA} for KASH and 
m-Algae~\cite{BrmA} for MAGMA, both by Nils Bruin; the TECC~\cite{Ki} calculator of Kida, 
also for KASH; 
and Simon's package~\cite{Sim} for Pari-GP. They are extremely useful in computing
ranks, and generators for the group $E(K)/2E(K)$. In each case below, it turns
out that the points generating $E(K)$ modulo $2 E(K)$ are actually generators
for the group $E(K)$ itself. This was proved using detailed height
calculations over the appropriate number field, with careful estimates for the
difference $\hat{h}(Q) - \frac{1}{2} h(Q)$ where $\hat{h}(Q)$ is the canonical 
height of the point $Q$, and $h(Q)$ the logarithmic height. The standard Silverman 
bounds~\cite{Sil2} are numerically too crude for our purposes, so recourse was made 
to the refinements of Siksek~\cite{Sik}.  Full details of the argument are given in 
an appendix to this paper~\cite{BTapp}. \\  
For the curve (\ref{E_1}) under immediate consideration, the cases 
$\delta=\eta_1,\eta_1 \eta_2$, give rise to curves of rank 0, and $\delta=1$, $\eta_2$, to 
curves of rank 1.\\
First, the curve (\ref{E_1}) at $\delta=1$ is
\begin{equation}
\label{ell1}
\mathcal{E}_1: Y^2 = X( X^2 -(\theta+\theta^2) X +(1+\theta+\theta^3))
\end{equation}
possessing only $2$-torsion over $K_1$, and with generator
\begin{equation}
\label{gen1}
G_1=(\frac{3+4\theta+\theta^2}{2},\;\frac{-4-6\theta-\theta^2-5\theta^3}{2}).
\end{equation}
The condition on the $X$-coordinate is
\[ -\theta+(3-3\theta+\theta^2-\theta^3) X \in \Q. \]
The point (\ref{gen1}) returns $(a,b)=(1,3)$.\\ \\ 
Second, the curve at $\delta=\eta_2$ is
\begin{equation}
\label{ell2}
\mathcal{E}_2: Y^2 = X( X^2 -(\theta-\theta^2) X +(1-\theta-\theta^3))
\end{equation}
possessing only $2$-torsion over $K_1$, with generator
\begin{equation}
\label{gen2}
G_2 = (\frac{1-\theta^2}{2},\;\frac{1-\theta}{2}).
\end{equation}
The condition on the $X$-coordinate is
\[ -\theta+(3+3\theta+\theta^2+\theta^3) X \in \Q. \]
The point (\ref{gen2}) returns $(a,b)=(1,1)$.\\ \\
Both curves (\ref{ell1}) and (\ref{ell2}) are minimal models.
\subsubsection{Equation (\ref{eq2}) and curves ${\mathcal E}_3,{\mathcal E}_4$} \label{curve2}
As above, (\ref{eq2}) leads to an equation
\[ (b+\frac{1}{\theta}a^2) (b^2+(-2+\frac{1}{\theta^2})a^4) = 
                                \pm \eta_1^i \eta_2^j (1+\theta) \Box, \]
where $i,j = 0,1$.  Specializing at the positive real root $0.643594...$ of $f_1(x)$, the sign 
of the right hand side must be positive. Putting $b/a^2 = \delta^{-1} x/(1+\theta)$, where 
$\delta=\eta_1^i \eta_2^j$, we thus have to find all $K_1$-points $(x,y)$ on the curves
\[(x+\frac{1+\theta}{\theta} \delta)(x^2+(-2+\frac{1}{\theta^2})(1+\theta)^2 \delta^2) 
                                                                              = y^2, \]
such that $\delta^{-1} \frac{x}{1+\theta} \in \Q$, for
$\delta=1$, $\eta_1$, $\eta_2$, or $\eta_1 \eta_2$.\\
Now put
\[ (x,y) = (2 X - \frac{1+\theta}{\theta} \delta,\; 2(1+\theta^2)Y), \]
to give
\begin{equation}
\label{E_2}
Y^2 = X(X^2 +(-1-2\theta-\theta^3) \delta X +(1+\theta+\theta^3) \delta^2);
\end{equation}
and the condition on the $X$-coordinate is
\[ \frac{2}{1+\theta} \frac{X}{\delta} - \frac{1}{\theta} \in \Q. \]
The cases $\delta=1$, $\eta_2$ give curves of rank 0; the remaining two cases are of rank 1.
First, the curve (\ref{E_2}) at $\delta=\eta_1$ is
\begin{equation}
\label{ell3}
\mathcal{E}_3: Y^2 = X(X^2 + (-1-\theta)X + (\theta+\theta^2-\theta^3)),
\end{equation}
possessing only $2$-torsion over $K_1$, with generator
\begin{equation}
\label{gen3}
G_3=(\frac{1-\theta^2}{2}, \; \frac{\theta^2+\theta^3}{2}).
\end{equation}
The condition on the $X$-coordinate is
\[  -2\theta-\theta^3 + (-3+7\theta-\theta^2+3\theta^3) X \in \Q. \]
The point (\ref{gen3}) returns $(a,b)=(1,1)$.\\ \\
Second, the curve (\ref{E_2}) at $\delta=\eta_1 \eta_2$ is
\begin{equation}
\label{ell4}
\mathcal{E}_4: Y^2 = X(X^2 + (-1+\theta)X + (-\theta+\theta^2+\theta^3)),
\end{equation} 
a conjugate of the curve (\ref{ell3}) under $\theta \rightarrow -\theta$.
Its generator is therefore
\begin{equation}
\label{gen4}
G_4 = (\frac{1-\theta^2}{2}, \; \frac{\theta^2-\theta^3}{2}).
\end{equation}
The condition on the $X$-coordinate is
\[ -2\theta-\theta^3 + (3+7\theta+\theta^2+3\theta^3) X \in \Q. \]
The point (\ref{gen4}) again returns $(a,b)=(1,1)$.
\subsubsection{Equation (\ref{eq3}) and curves ${\mathcal E}_i, i=5,\ldots,8 $} \label{curve3}
As above, (\ref{eq3}) leads to an equation of type
\[ (b+2\phi a^2) (b^2 +4 (4+\phi^2)a^4) = \pm \epsilon_1^i \epsilon_2^j \Box, \]
for $i,j=0,1$.  Specializing $\phi$ at the positive real root 
$0.910179...$ of $f_2(x)$,
it follows that the sign must be positive. For $a \neq 0$, put $b/a^2 = \delta^{-1} x$,
where $\delta=\epsilon_1^i \epsilon_2^j$, which leads to
seeking all $K_2$-points $(x,y)$ on the curves
\[ (x+2 \phi \delta) (x^2 +4(4+\phi^2) \delta^2) = y^2, \]
subject to $\delta^{-1}x \in \Q$, with $\delta = 1$, $\epsilon_1$, $\epsilon_2$, or
$\epsilon_1 \epsilon_2$, that is, $\delta=1$, $\frac{1}{2}\phi+\frac{1}{4}\phi^3$, 
$2+2\phi+\frac{1}{2}\phi^2+\frac{1}{2}\phi^3$, 
            $1+\frac{3}{2}\phi+\frac{1}{2}\phi^2+\frac{1}{4}\phi^3$. 
Put 
\[ (x,y) = (4 X-2\phi\delta,\; (2+\phi^2)^2 Y) \]
to give
\begin{equation}
\label{E_3}
Y^2 = X(X^2 -\phi \delta X +(1+\frac{1}{2}\phi^2)\delta^2),
\end{equation}
where the condition on $X$-coordinate has become 
\[ -2\phi+\frac{4}{\delta}X \in \Q. \]
All four curves are of rank 1.
The curve (\ref{E_3}) with $\delta=1$ has equation
\begin{equation}
\label{ell5}
\mathcal{E}_5: Y^2 = X(X^2 -\phi X +(1+\frac{1}{2}\phi^2)),
\end{equation}
possessing only $2$-torsion over $K_2$, with generator
\begin{equation}
\label{gen5}
G_5 = (2-2\phi+\frac{1}{2}\phi^2-\frac{1}{2}\phi^3, \; 5-5\phi+\phi^2-\phi^3).
\end{equation}
The condition on the $X$-coordinate is 
\[ -2\phi+4 X \in \Q. \]
Twice the generator at (\ref{gen5}) is the point
\[ (X,Y) = (\frac{5}{4}+\frac{1}{2}\phi, 
         \frac{1}{2}+\frac{7}{4}\phi-\frac{3}{16}\phi^2+\frac{5}{16}\phi^3), \]
which leads to $(a,b)=(1,5)$.\\ \\
The curve (\ref{E_3}) with $\delta=\epsilon_1$ has equation
\begin{equation}
\label{ell6}
\mathcal{E}_6: Y^2=X(X^2 +(-1+\frac{1}{2}\phi^2)X+(1-\frac{1}{2}\phi^2)),
\end{equation}
possessing only $2$-torsion over $K_2$, with generator
\begin{equation}
\label{gen6}
G_6 = (1-\frac{1}{2}\phi^2, \; 1-\frac{1}{2}\phi^2).
\end{equation}
The condition on the $X$-coordinate is
\[ -2\phi+(6+\phi^3) X \in \Q. \]
The curve (\ref{E_3}) with $\delta=\epsilon_2$ has equation
\begin{equation}
\label{ell7}
\mathcal{E}_7: 
  Y^2 = X(X^2 +(-2-2\phi-\frac{1}{2}\phi^3) X+(13+14\phi+\frac{5}{2}\phi^2+3\phi^3)),
\end{equation}
possessing only $2$-torsion over $K_2$, with generator
\begin{equation}
\label{gen7}
G_7 = (1+\frac{1}{2}\phi+\frac{1}{4}\phi^3, \; -3-3\phi-\frac{1}{2}\phi^2-\frac{1}{2}\phi^3).
\end{equation}
The condition on $X$-coordinate has become
\[ -2\phi+2(4-4\phi+\phi^2-\phi^3) X \in \Q. \]
The curve (\ref{E_3}) with $\delta=\epsilon_1 \epsilon_2$ has equation
\begin{equation}
\label{ell8}
\mathcal{E}_8: Y^2 = 
X(X^2 +(-1-\phi-\frac{1}{2}\phi^2-\frac{1}{2}\phi^3)X+(5+6\phi+\frac{3}{2}\phi^2+\phi^3)),
\end{equation}
possessing only $2$-torsion over $K_2$, with generator
\begin{equation}
\label{gen8}
G_8 = (1+\frac{1}{2}\phi+\frac{1}{4}\phi^3, \; -2-2\phi-\frac{1}{2}\phi^3).
\end{equation}
The condition on $X$-coordinate
has become 
\[ -2\phi +(-12+14\phi-2\phi^2+3\phi^3) X \in \Q. \]
All curves are minimal models.
\subsubsection{Equation (\ref{eq4}) and curves ${\mathcal E}_i, i=9,\ldots, 12$} \label{curve4}
As in the third case, we deduce an equation in $\mathcal{O}_2$:
\[ (b+\frac{4}{\phi} a^2) (b^2 +(-16+\frac{16}{\phi^2}) a^4) = 
                                      \pm \epsilon_1^i \epsilon_2^j \Box, \]
where $i,j = 0,1$, and specializing at the positive real root of $f_2(x)$, the sign must be 
positive.
For $a \neq 0$, put $b/a^2 = \delta^{-1} x$, where $\delta=\epsilon_1^i \epsilon_2^j$.
This leads to finding all $K_2$-points $(x,y)$ on the curves 
\[ (x+\frac{4}{\phi} \delta) (x^2 +(-16+\frac{16}{\phi^2}) \delta^2) = \Box, \]
with
$\delta^{-1} x \in \Q$, and $\delta=1$, $\epsilon_1$, $\epsilon_2$, $\epsilon_1 \epsilon_2$.
Putting 
\[ (x,y)=(4X -\frac{4}{\phi} \delta, \; (2+\phi^2)^2 Y ) \]
gives
\begin{equation}
\label{E_4}
Y^2 = X(X^2 -\frac{2}{\phi} \delta X + (-1+\frac{2}{\phi^2}) \delta^2 ),
\end{equation}
and the condition on $X$ becomes
\[ -\frac{4}{\phi} +\frac{4}{\delta} X \in \Q. \]
The curve (\ref{E_4}) with $\delta=1$ has equation
\begin{equation}
\label{ell9}
\mathcal{E}_9: Y^2=X(X^2 +(-2\phi-\frac{1}{2}\phi^3)X+(1+\frac{1}{2}\phi^2)),
\end{equation}
of rank 1, 
possessing only $2$-torsion over $K_2$, with generator
\begin{equation}
\label{gen9}
G_9 = (1+\frac{1}{2}\phi+\frac{1}{4}\phi^3, -\phi).
\end{equation}
of canonical height $0.125726743336419...$ 
The condition on $X$ has become
\[ -4\phi-\phi^3 + 4 X \in \Q. \]
The curve (\ref{E_4}) with $\delta=\epsilon_1$ has equation
\begin{equation}
\label{ell10}
\mathcal{E}_{10}: Y^2 = X(X^2 +(-1-\frac{1}{2}\phi^2)X+(1-\frac{1}{2}\phi^2)),
\end{equation}
and is of rank 2, possessing only $2$-torsion over $K_2$, with generators 
\begin{equation}
\label{gen101}
P_1: (X,Y) = (1, \frac{1}{2}\phi^2),
\end{equation}
and 
\begin{equation}
\label{gen102}
P_2: (X,Y) = (\frac{1}{2}\phi+\frac{1}{2}\phi^2-\frac{1}{4}\phi^3, 1-\frac{3}{2}\phi^2).
\end{equation}
The condition on $X$ is 
\begin{equation} \label{X-condition 4.2}
-4\phi-\phi^3+(6\phi+\phi^3)X \in \Q
\end{equation}
The curve (\ref{E_4}) with $\delta=\epsilon_2$ has equation
\begin{equation}
\label{ell11}
\mathcal{E}_{11}: Y^2 = 
                X(X^2 +(-4-5\phi-\phi^2-\phi^3)X+(13+14\phi+\frac{5}{2}\phi^2+3\phi^3)),
\end{equation}
of rank 1, possessing only $2$-torsion over $K_2$, with generator
\begin{equation}
\label{gen11}
G_{11} = (2+2\phi+\frac{1}{2}\phi^2+\frac{1}{2}\phi^3, -2-2\phi-\frac{1}{2}\phi^2-\frac{1}{2}\phi^3).
\end{equation}
The condition on $X$ becomes 
\[ -4\phi-\phi^3+(8-8\phi+2\phi^2-2\phi^3) X \in \Q. \]
The curve (\ref{E_4}) with $\delta=\epsilon_1 \epsilon_2$ has equation
\begin{equation}
\label{ell12}
\mathcal{E}_{12}: Y^2 = 
X( X^2 +(-3-3\phi-\frac{1}{2}\phi^2-\frac{1}{2}\phi^3)X+(5+6\phi+\frac{3}{2}\phi^2+\phi^3)),
\end{equation}
of rank 1, 
possessing only $2$-torsion over $K_2$, with generator
\begin{equation}
\label{gen12}
G_{12} = (1+\frac{1}{2}\phi+\frac{1}{4}\phi^3, -1-\phi-\frac{1}{2}\phi^2-\frac{1}{2}\phi^3).
\end{equation}
The condition on $X$ has become
\[ -4\phi-\phi^3+(-12+14\phi-2\phi^2+3\phi^3) X \in \Q. \]
All curves are minimal models.\\ \\ \\

\newcommand{\ser}[1]{\langle{#1}\rangle} 
\newcommand{\serd}[2]{\langle{#1,#2}\rangle} 

\newtheorem{Results_r1}{Proposition}
\newtheorem{fact2}[Results_r1]{Fact}
\newtheorem{a la Skolem}[Results_r1]{Theorem}
\newtheorem{rational_mod}[Results_r1]{Proposition-Definition}
\newtheorem{algorithm2}[Results_r1]{Algorithm}
\newtheorem{Results_r2}[Results_r1]{Proposition}

\subsection{Cases corresponding to rank 1 elliptic curves}\label{solutions_r1}
We gave a detailed discussion of the solution of Problem \ref{problem} for rank one
elliptic curves in section 4 of our companion paper \cite{BT}, in which we also gave a 
number of concrete examples. Therefore, we confine ourselves here in giving all necessary 
data for the corresponding rank one elliptic curves of section \ref{elliptic curves} and 
saying that, following exactly the same method and working $p$-adically with $p=3$, we conclude 
the following result:    
\begin{Results_r1} \label{Results_r1}
For each elliptic curve ${\mathcal E}_i,\,i=1,\ldots,9$, and $i=11,12$, the only points on 
${\mathcal E}_i$ whose $X$-coordinate belongs to the appropriate quartic field and which satisfies 
the corresponding condition $\beta X+\gamma \in\Q$, are given by the following:
\begin{itemize}
\item Elliptic curve ${\mathcal E}_1: \mbox{points}\; \pm G_1$, giving $a=\pm 1,b=3$ at 
(\ref{eq1}). \\ From section \ref{odd P}, $P=1,Q=-4$.
\item Elliptic curve ${\mathcal E}_2: \mbox{points}\; \pm G_2$, giving $a=\pm 1,b=1$ at 
(\ref{eq1}). \\ From section \ref{odd P}, $P=1,Q=1$.
\item Elliptic curve ${\mathcal E}_3: \mbox{points}\; \pm G_3$, giving $a=\pm 1,b=-1$ at
(\ref{eq2}). \\ From section \ref{odd P}, $P=1,Q=1$.
\item Elliptic curve ${\mathcal E}_4: \mbox{points}\; \pm G_4$, giving $a=\pm 1,b=1$ at
(\ref{eq2}). \\ From section \ref{odd P}, $P=1,Q=1$.
\item Elliptic curve ${\mathcal E}_5: \mbox{points}\; \pm 2G_5$, giving $a=\pm 1,b=5$ at
(\ref{eq3}). \\ From section \ref{even P}, $P=4,Q=-17$.
\item Elliptic curve ${\mathcal E}_6: \mbox{no point}$.
\item Elliptic curve ${\mathcal E}_7: \mbox{points}\; \pm 2G_7$, giving $a=\pm 1,b=2$ at
(\ref{eq3}). \\ From section \ref{even P}, $P=4,Q=4$, 
rejected (we assumed $P,Q$  relatively prime).  
\item Elliptic curve ${\mathcal E}_8: \mbox{points}\; \pm 2G_8$, giving $a=\pm 1,b=0$ at
(\ref{eq3}), which is impossible.
\item Elliptic curves ${\mathcal E}_9, {\mathcal E}_{11}, {\mathcal E}_{12}: 
\mbox{no points}$. 
\end{itemize}
\end{Results_r1}
\subsection{Cases corresponding to rank 2 elliptic curves}\label{solutions_r2}
For the solution of Problem \ref{problem} when the rank of the elliptic curve is 2,
we make the following assumptions:

\vspace{1ex}
{\em Assumption 1}. There exists a rational prime $p$ with the following properties: 
\vspace{-2mm}
\begin{itemize}
\item $f(X)$ is irreducible in $\Q_p[X]$. This implies that $p$ is a prime divisor of the 
number
field $\Q(\alpha)$ and there is only one discrete 
(normalized) valuation $v$ defined on $\Q(\alpha)$ with $v(p)=1$. Moreover, the completion of 
$\Q(\alpha)$ with respect to $v$ is $\Q_p(\alpha)$ and, 
according to our assumptions, $[\Q_p(\alpha):\Q_p]=[\Q(\alpha):\Q]=d $.
\item The coefficients of (\ref{gen_xy_elliptic}) are in $\Z_p[\alpha]$.
\item Equation (\ref{gen_xy_elliptic}) is a minimal Weierstrass equation for 
$\cE/\Q_p(\alpha)$ at $v$. 
\item $\beta,\gamma \in\Q_p(\alpha)$ are $p$-adic units. 
\end{itemize}

\vspace{1ex}
{\em Assumption 2}. We know two independent points $Q_1,Q_2\in\cE(\Q(\alpha))$, each having the 
form $(s/t^2,u/t^3)$ with $s,u\in\Z[\alpha]$, $t$ a positive integer divisible by $p$ and
$(\Nm(s),t)=(\Nm(u),t)=1$; here $\Nm$ denotes norm relative to the extension $\Q(\alpha)/\Q$. 
If $p=2$ we assume something more, namely, that $t$ is divisible by $p^2=4$.  

According to the notation and facts in section 4 of our paper \cite{BT}, 
$Q_i\in\hat{\cE}({\cal M}^r),\, (i=1,2)$. The same arguments used therein, lead to the 
following conclusion:
\begin{fact2} \label{fact2}
Let $P=(X_0,Y_0)$ be any finite point of $\cE(\Q)$ and let $n_1,n_2$ denote integer variables. 
Then, both $\beta x(P+n_1Q_1+n_2Q_2)+\gamma$ and 
$\displaystyle{(\beta x(n_1Q_1+n_2Q_2)+\gamma)^{-1}}$ can be expressed as   
$\theta_0(n_1,n_2)+\theta_1(n_1,n_2)\alpha+\cdots +\theta_{d-1}(n_1,n_2)\alpha^{d-1}$,  
where each $\theta_i(n_1,n_2)$ is a $p$-adically convergent power series in $n_1,n_2$ with
coefficients in $\Z_p$, having also the following property: For every $(k,\ell)\neq (0,0)$,
\begin{equation} \label{v_of_coeff_2}
v(\mbox{coefficient of \,}n_1^kn_2^{\ell})\geq 
\begin{cases} 
\left\lfloor\frac{(p-2)(k+\ell)}{p-1}\right\rfloor +1 & \mbox{if $p\geq 3$} \\
                                                          k+\ell+1 & \mbox{if $p=2$}  
 \end{cases} \;.
\end{equation}
\end{fact2}

The coefficients of the series $\theta_i$ depend on the coordinates of $Q_1,Q_2$ and, in case 
of $\beta x(P+n_1Q_1+n_2Q_2)+\gamma$, also on the coordinates of $P$.

\vspace{1ex}
{\em Assumption 3}. The typical point on $\cE(\Q(\alpha))$ can be expressed in the form 
$P+n_1Q_1+n_2Q_2$, where $P$ is chosen from a finite explicitly known set of points, including 
the zero point. 

\vspace{1ex}
Under Assumptions 1-3, problem \ref{problem} is clearly reduced to solving the system of 
equations 
$\theta_1(n_1,n_2)=0,\ldots,\theta_{d-1}(n_1,n_2)=0$ for each value of $P$. In \cite{BT} 
we had a similar problem, but for a curve of rank 1, and the system of equations we had 
to solve was in one unknown $n_1$. In that situation, Strassman's theorem (see, for example, 
Theorem 4.1 in \cite{BT}) was applicable, but not in the present one, where we have two 
unknowns $n_1,n_2$. Instead, we apply a
theorem, which we state and prove below, inspired by the paper of Th.~Skolem \cite{Sk}. 

It is worth mentioning that, in a similar situation, S.~Duquesne in \cite{D} applied a 
different method based on his explicit version of a $p$-adic Weierstrass preparation 
theorem of T.~Sugatani \cite{Su} (see sections 2 and 3 of \cite{D}). That explicit version of 
Sugatani's theorem is interesting, but from our experience (in a first unpublished 
version of this paper, we employed Duquesne's method) its application is more complicated. 
  
Our remarks a few lines above make evident that, in order to solve problem \ref{problem}, we 
must know how to find explicitly all $p$-adic integer solutions of a system of equations 
$F_1=0,F_2=0$, for appropriate series $F_1,F_2\in\Z_p[[x_1,x_2]]$.   
In a more general setting we state and prove the theorem below which we will apply in the 
special case of two unknowns.
\begin{a la Skolem} \label{a la Skolem}
Let $p$ be a prime and for $r=1,\ldots,n$ let 
\[F_r(x_1,\ldots,x_n)=\sum_{i=0}^{\infty} p^i f_{ir}(x_1,\ldots,x_n)\:,\]
where $f_{ir}\in\Z_p[x_1,\ldots,x_n]$. Assume that the following conditions are satisfied:
\begin{enumerate}
\item $f_{0r}(x_1,\ldots,x_n)$ is homogeneous of degree, say, $d_r\geq 1$.
\item Every monomial $x_1^{i_1}\cdots x_n^{i_n}$ in $F_r(x_1,\ldots,x_n)$ is of degree at least 
$d_r$ (this, in particular, implies $F_r(0,\ldots,0)=0$).
\item For every $r=1,\ldots,n$ there exist $h_{1r}, \ldots, h_{nr}\in\Z_p[x_1,\ldots,x_n]$ such 
that $h_{1r}\cdot f_{01}+\cdots +h_{nr}\cdot f_{0n}=H_r\in\Z_p[x_r]$
and the only solution to the congruence $H_r(x)\equiv 0\pmod{p}$ is $x\equiv 0\pmod{p}$ 
(this, in particular, implies that $H_r$ is a non-zero polynomial $\bmod{p}$). 
\end{enumerate}
Then, the only solution in $p$-adic integers of the system 
$F_r(x_1,\ldots,x_n)=0,\, (r=1,\ldots,n)$ is the zero solution. 
\end{a la Skolem}
\Proof 
Suppose $F_r(x_1,...,x_n)=0$ for $r=1,...,n$, where $x_i \in {\bf Z}_p$
are not all zero. Then
\[ f_{0r}(x_1,...,x_n) \equiv 0 \pmod p, \quad r=1,...,n \]
so that by hypothesis (3),
\[ H_r(x_r) \equiv 0 \pmod p, \quad r=1,...,n, \]
that is, also by hypothesis (3),
\[ x_r \equiv 0 \pmod p, \quad r=1,...,n. \]
Thus $p^{\alpha_r} || x_r$, $\alpha_r \geq 1$ (with convention that
$\alpha_r = \infty$ if $x_r = 0$). Define the integer $j$ in the range
$1 \leq j \leq n$ by $\alpha_j=\min{(\alpha_1,\alpha_2,...,\alpha_n)} = \alpha$.
(The integer $j$ exists since at least one $\alpha_r$ is finite).
Now put $x_r = p^\alpha X_r$, $r=1,...,n$, where $X_r \in {\bf Z}_p$, 
and $p \not | X_j$.\\ \\
Then
\[ F_r(x_1,...,x_n) = 0, \quad r=1,...,n, \]
implies
\[ \sum_{i=0}^\infty p^i f_{ir}(p^\alpha X_1,...,p^\alpha X_n) = 0, \quad r=1,...,n, \]
that is,
\[ p^{\alpha d_r} [ f_{0r}(X_1,...X_n) + \sum_{i=1}^\infty p^i g_{ir}(X_1,...,X_n) ] = 0, 
                                                                            \quad r=1,...,n, \]
where $g_{ir}(X_1,...,X_n) \in {\bf Z}_p[X_1,...,X_n]$, using hypotheses (1) and (2).
Thus 
\[ f_{0r}(X_1,...X_n) + \sum_{i=1}^\infty p^i g_{ir}(X_1,..,X_n)) = 0, \quad r=1,...,n, \]
so that
\[ f_{0r}(X_1,...,X_n) \equiv 0 \pmod p, \quad r=1,...,n, \]
whence by hypothesis (3), 
\[ H_r(X_r) \equiv 0 \pmod p, \quad r=1,...,n. \]
In particular,
\[ H_j(X_j) \equiv 0 \pmod p, \]
so that $X_j \equiv 0 \pmod p$, by hypothesis (3), contrary to assumption. \proofend \\[3mm]

{\bf Remarks} (1) If for every $r=1,\ldots,n$, $d_r=1$, hence 
$f_{0r}= a_{1r}x_1+\cdots+a_{nr}$ (say), then the conditions of the theorem are equivalent to 
the non-vanishing $\bmod{p}$ of the determinant of the matrix $(a_{ir})$. \\
(2) When $n=2$, at least the existence of the polynomials $h_{1r},h_{2r},\, (r=1,2)$ is 
guaranteed by the basic theory of resultants; in that case, $H_1(x_1)$ is the resultant of the 
polynomials
$f_{01}(x_1,x_2),f_{02}(x_1,x_2)$ with respect to the variable $x_2$, and analogously for 
$H_2(x_2)$.     
%
\subsubsection*{Application of Theorem \ref{a la Skolem} to (\ref{ell10})}
A Mordell-Weil basis for the elliptic curve (\ref{ell10}) over $\Q(\phi)$  is formed by the 
generators of infinite order 
$P_1=(1,\frac{1}{2}\phi^2), P_2=(\frac{1}{2}\phi+\frac{1}{2}\phi^2-\frac{1}{4}\phi^3, 
                                                                       1-\frac{3}{2}\phi^2)$ 
(see section 4.10 in the appendix to~\cite{BTapp}) and the generator
$T=(0,0)$ of the torsion subgroup.
We define $Q_1=P_1+8P_2$, 
$Q_2=24P_2$. Note that $\{Q_1, P_2\}$ remains a basis for the torsion-free part of the group of
rational points of (\ref{ell10}) over $\Q(\phi)$, therefore any non-zero point 
$(X,Y)\in {\mathcal E}_{10}(\Q(\phi))$ can be written as 
\begin{equation} \label{generic point}
kP_2+\epsilon T+n_1Q_1+n_2Q_2\,,\;\; n_1,n_2\in\Z\,, \; k\in\{-11,\ldots,12\}\,, \; 
                                                                    \epsilon\in\{0,1\} \;,
\end{equation}
and $n_1,n_2,k,\epsilon$ not all zero. 

Note that Assumption 1 at the beginning of section \ref{solutions_r2} is fulfilled with 
$p=3$ and $\beta=6\phi+\phi^3,\gamma=-4\phi-\phi^3$.
Assumption 2 is then fulfilled for the points $Q_1,Q_2$ defined above.  
In (\ref{generic point}) we put $P=kP_2+\epsilon T$. 
There are $24\cdot 2=48$ possibilities for $P$, with points other than for $k=\epsilon=0$ 
being ``finite points". The generic point $(X,Y)\in {\mathcal E}_{10}(\Q(\phi))$ 
has the form $P+n_1Q_1+n_2Q_2$, and hence Assumption 3 is also fulfilled.
We are interested in finding all points $(X,Y)$ as above, that satisfy condition (\ref{X-condition 4.2}).
Therefore, if at least one of $k,\epsilon$ is non-zero, we may assume, since $T = -T$,
that $k\in\{1,\ldots,12\}$ if $\epsilon=0$ 
and 
$k\in\{0,\ldots,12\}$ if $\epsilon=1$, reducing thus to $1+12+13=26$ the possibilities for 
the point $P$.  

Following the method of Flynn and Wetherell \cite{FW} as described in section 4 of \cite{BT}, 
we have (in the notation of \cite{BT})
\begin{eqnarray*}
z(Q_1) & = &  33+240\phi+33\phi^2+93\phi^3 +O(3^5) \in {\cal M}\\
z(Q_2) & = &  213+234\phi+105\phi^2+144\phi^3+O(3^5) \in {\cal M}\;.
\end{eqnarray*}
The ``addition law" in the {\em formal group} of our elliptic curve is given by
\[
\begin{aligned}
 \cF(z_1,z_2) = & z_1+z_2 + (\frac{1}{2}\phi^2+1)z_1z_2^2+(\frac{1}{2}\phi^2+1)z_1^2z_2
    + (\phi^2-2)z_1z_2^4 \\
    &+(2\phi^2-2)z_1^2z_2^3+(2\phi^2-2)z_1^3z_2^2+(\phi^2-2)z_1^4z_2+\cdots
\end{aligned}
\]
The logarithmic and exponential series in the formal group are
\[
\begin{aligned}
\log t  = &\,  t+(-\frac{1}{6}\phi^2-\frac{1}{3})t^3+\frac{2}{5}t^5+O(t^7)\\
 \exp t = &\, t+(\frac{1}{6}\phi^2+\frac{1}{3})t^3+\frac{4}{15}t^5+O(t^7)
\end{aligned}
\]
For any point $Q$ on the elliptic curve we will use the notation $X(Q)$ for the $X$-coordinate 
of the point $Q$.  
For any finite points $P=(X_0,Y_0)$ and $R$ of our elliptic curve we express 
$\beta X(P+R)+\gamma$
(with $\beta=6\phi+\phi^3$ and $\gamma=-4\phi-\phi^3$)
as a formal power 
series of $z(R)$ with coefficients in $\Z[\phi,X_0,Y_0]$: 
{\small
\begin{equation}\label{bx+g}
\begin{aligned} 
\beta X(P+R) +\gamma = & \,(X_0-1)\phi^3+(6X_0-4)\phi +(2Y_0\phi^3+12Y_0\phi)\,z(R)\\
    & +[(3X_0^2-4X_0)\phi^3+(4-16X_0+18X_0^2)\phi]\,z(R)^2 \\
       & + [(4Y_0X_0-4Y_0)\phi^3+(24Y_0X_0-16Y_0)\phi]\,z(R)^3 \\
    & + [(4X_0-2+Y_0^2+4X_0^3-12X_0^2)\phi^3 \\
  & \hspace{6mm}+(24X_0^3-48X_0^2+32X_0-4+6Y_0^2)\phi]\,z(R)^4
 + O(z(R)^5)
 \end{aligned}
\end{equation}  
}
We also express the inverse of $\beta X(R)+\gamma$ as a formal power series in $z(R)$:
{\small
\begin{equation} \label{inverse_bx+g}
\frac{1}{\beta X(R)+\gamma} =  \frac{\phi^3+2\phi}{16}z(R)^2
    -\frac{\phi}{8}z(R)^4+\frac{5\phi^3+2\phi}{32}z(R)^6 +O(z(R)^8)
\end{equation}
}
We have the 3-adic expansions
\[
\begin{aligned}
\log z(Q_1) = &\, 3(32+35\phi+50\phi^2+61\phi^3)+O(3^5) \in 3\Z_3[\phi]\\
\log z(Q_2) = &\, 3(47+38\phi^2)+3^2(8\phi+7\phi^3)+O(3^5) \in 3\Z_3[\phi]\;.
\end{aligned}
\]
Let $n_1,n_2$ be integers and set $R=n_1Q_1+n_2Q_2$. From section 4 of \cite{BT} we 
know that 
\vspace{-3mm}
\[z(R)=z(n_1Q_1+n_2Q_2)=\exp(n_1\log z(Q_1)+n_2\log z(Q_2))\in 3\Z_3\serd{n_1}{n_2}[\phi] \;.\]  
This can be easily computed $\bmod{\,3^5}$; we need consider only the first three terms of the
exponential series, in view of the fact that 
$\log z(Q_1),\log z(Q_2)\in 3\Z_3[\phi]$.
{\small    
\begin{equation} \label{z_coord}
\begin{aligned}
z(R)\bmod{\,3^5}=&\, 216n_1n_2^2+81n_1^2n_2+96n_1+141n_2+180n_1^3+153n_2^3+81n_2^4n_1\\
   & \hspace{1cm}+162n_2^5+81n_2n_1^4+162n_2^2n_1^3 
   \\
&\,+(135n_1^3+162n_2^3+72n_2+105n_1+81n_2^4n_1+216n_1n_2^2\\
   & \hspace{1cm}+81n_2n_1^4+81n_2^3n_1^2+108n_1^2n_2)\phi \\
&\,+(150n_1+114n_2+72n_2^3+162n_2^5+81n_2^4n_1+135n_1n_2^2\\
 & \hspace{1cm} +126n_1^3+108n_1^2n_2+81n_2^3n_1^2+162n_2^2n_1^3)\phi^2 \\
&\,+(81n_2^3+81n_1^5+183n_1+63n_2+72n_1^3+162n_2^3n_1^2\\
 & \hspace{1cm} +162n_2^4n_1+162n_2^2n_1^3+135n_1^2n_2+189n_1n_2^2)\phi^3
\end{aligned}
\end{equation}
}
As noted in (\ref{fact2}), substitution of the above value for $z(R)$ in (\ref{bx+g}) and (\ref{inverse_bx+g}) gives, 
after reduction $\bmod{\,3^5}$, an element in $\Z\serd{n_1}{n_2}[\phi,X_0,Y_0]$ and 
$\Z\serd{n_1}{n_2}[\phi]$, respectively (the formulas are too 
long, especially the first one, to be included here). This is of the form
\begin{equation} \label{theta expression} 
 \theta_0(n_1,n_2)+ \theta_1(n_1,n_2)\phi + \theta_2(n_1,n_2)\phi^2 + 
                                                 \theta_3(n_1,n_2)\phi^3 \;,
\end{equation}
where $\theta_i(x,y)\in\Z[x,y]$ and, in the first case, with coefficients depending on 
$X_0,Y_0$. 
\begin{quote}
{\footnotesize
{\em Notation}. In the sequel we assume that $(X,Y)$ is a point on the curve 
${\mathcal E}_{10}$, 
such that $X$ satisfies condition (\ref{X-condition 4.2}). We put $R=n_1Q_1+n_2Q_2$, with 
$n_1,n_2\in\Z$. Note that, the typical form of $(X,Y)$ is either $(X,Y)=P+R$ with $P=(X_0,Y_0)$
belonging to the set of 25 ``finite" points mentioned at the beginning of this section, or 
$(X,Y)=R$. 
}
\end{quote}

\vn
\underline{Case 1}: $\;(X,Y)=P+R$.
We recall that $P=(X_0,Y_0)=kP_2+\epsilon T,\, k=0,1,\ldots,12\,,\epsilon=0,1$. 
Suppose first $\epsilon=1$. Using the computer we find, for every specific $P$, an explicit 
expression for the form (\ref{theta expression}) 
for $\beta X +\gamma \pmod{3^5}$. 
In every case but $k=4$, we find out that $\theta_i(n_1,n_2)\not\equiv 0\pmod{3}$ for at 
least one 
$i$, hence $\beta X +\gamma$ cannot be a rational number. When $k=4$, we compute 
$\theta_1(n_1,n_2)\equiv 6+6n_1+6n_2\pmod{3^2}$
and 
$\theta_3(n_1,n_2)\equiv 3n_1+3n_2\pmod{3^2}$, therefore the simultaneous vanishing of 
$\theta_1(n_1,n_2)$ and $\theta_3(n_1,n_2)$ is impossible. 
This leads to the conclusion that $\beta X+\gamma$ cannot be a rational number.
Next, consider the case $\epsilon=0$. In every case but $k=2,10$, we see that 
$\theta_i(n_1,n_2)\not\equiv 0\pmod{3}$ for at least one $i$, hence $\beta X +\gamma$ 
cannot be a rational number. 

The cases $k=2,10$ need a deeper treatment. Working $p$-adically with $p=3$ we apply 
Theorem \ref{a la Skolem} in order to solve in 3-adic integers the system 
\begin{equation} \label{system1}
\theta_3(n_1,n_2)=0,\: \theta_2(n_1,n_2)=0 \;\; n_1,n_2 \in\Z_3 \;.
\end{equation}

\underline{Case 1.1}: $P=2P_2$. 
We are looking for points $(X,Y)=2P_2+n_1Q_1+n_2Q_2$ such that $X$ satisfies 
condition (\ref{X-condition 4.2}). 
Note that, for $(n_1,n_2)=(0,0)$ this is satisfied. Indeed, then
\[(X,Y)=2P_2=(\frac{1}{2}-\frac{1}{2}\phi+\frac{1}{4}\phi^2-\frac{1}{4}\phi^3, 
                                      \frac{1}{2}-\frac{1}{4}\phi-\frac{1}{8}\phi^3) \]
and we check that $\beta X+\gamma=-4$, as required. This means that $(n_1,n_2)=(0,0)$ is a 
solution
to the system (\ref{system1}). 
Keeping in mind this solution we define 
\[F_1(n_1,n_2)=\frac{1}{3}\theta_3(n_1,n_2)\,,\: F_2(n_1,n_2)=\frac{1}{3}\theta_2(n_1,n_2)\]
and, using theorem \ref{a la Skolem}, we will show that $(n_1,n_2)=(0,0)$ is the only solution 
of the system $F_1=0,F_2=0$ in 3-adic integers. We compute 
\begin{align*}
F_1(n_1,n_2) = &  2n_1 + 3(n_1^3+n_1^2+n_1+n_2+2n_2^2)\\
   & + 3^2(n_2+2n_1+n_2^3+2n_2^2+n_2^4+2n_1n_2+2n_1n_2^2) +3^3(\cdot) \;, 
\end{align*}
where $(\cdot)$ denotes a series in $\Z\langle n_1,n_2\rangle$ with zero constant term. Also,
\begin{align*}
F_2(n_1,n_2) = & n_1+n_2 +3(2n_1^3+n_1^2+2n_1n_2+n_1+n_2^3) \\
    &+ 3^2(2n_1^2+2n_2^3+2n_2^2+2n_1^2n_2+2n_1n_2+n_1n_2^2+2n_1n_2^3) +3^3(\cdot)\;,
\end{align*}
where $(\cdot)$ is as above. Actually the essential terms are 
$f_{01}=2n_1$ and $f_{02}=n_1+n_2$, 
with corresponding determinant of their coefficients
\[ \left| \begin{array}{cc} 2 & 0 \\ 1 & 1 \end{array} \right| \;. \]
This is non-zero $\bmod{\,3}$, hence, by remark (1) following theorem \ref{a la Skolem}, 
the only solution to our system is $(n_1,n_2)=(0,0)$.
This corresponds to the point $2P_2$ on the curve ${\mathcal E}_{10}$ with $X$-coordinate 
$\frac{1}{2}-\frac{1}{2}\phi+\frac{1}{4}\phi^2-\frac{1}{4}\phi^3$. 
Then, in section \ref{even P} $(a,b)=(1,-4)$ which does not furnish us with a solution of 
equation (\ref{eq4}). 

\vspace{1mm}
\underline{Case 1.2}: $P=10P_2$. 
Now we are looking for points $(X,Y)=10P_2+n_1Q_1+n_2Q_2$ such that $X$ satisfies 
condition (\ref{X-condition 4.2}). 
Note that, for $(n_1,n_2)=(2,-1)$ the condition is satisfied. Indeed, then
\[(X,Y)=10P_2+2Q_1-Q_2=2P_1+2P_2=(\frac{1}{2}+\frac{1}{2}\phi+\frac{1}{4}\phi^2
                     +\frac{1}{4}\phi^3, -\frac{1}{2}-\frac{1}{4}\phi-\frac{1}{8}\phi^3) \]
and we check that $\beta X+\gamma=4$, as required. In particular, we conclude that 
$(n_1,n_2)=(2,-1)$ is a solution to (\ref{system1}). 
Therefore, we put $n_1=x_1+2,n_2=x_2-1$, we define  
\begin{align*}
F_1(x_1,x_2)= & \frac{1}{3}\theta_3(n_1,n_2) =\frac{1}{3}\theta_3(x_1+2,x_2-1)\,,\\
F_2(x_1,x_2)= &\frac{1}{3}\theta_2(n_1,n_2) =\frac{1}{3}\theta_2(x_1+2,x_2-1)
\end{align*}
and we will show, using theorem \ref{a la Skolem}, that $(x_1,x_2)=(0,0)$ is the only solution 
in 3-adic integers to the system $F_1=0,F_2=0$. We compute 
\begin{align*} 
F_1(x_1,x_2) = & 2x_1+3(2x_2x_1+2x_2^2+x_1^3)+ 3^2(\cdot) \\
F_2(x_1,x_2) = & x_1+x_2 +3(2x_2x_1+x_2^3) +3^2(\cdot)\;,      
\end{align*}
where $(\cdot)$ denotes a series in $\Z\langle n_1,n_2\rangle$ with zero constant term.
As in case 1.1, the determinant of the coefficients of the first-degree terms $2x_1$ and 
$x_1+x_2$ is non-zero $\bmod{3}$, therefore $(x_1,x_2)=(0,0)$ is the only solution of the 
system in 3-adic integers. It follows that, in case 1.2, $(n_1,n_2)=(2,-1)$ is the only 
possible solution of the system (\ref{system1}).  
This gives a point on the curve ${\mathcal E}_{10}$ with $X$-coordinate 
$\frac{1}{2}+\frac{1}{2}\phi+\frac{1}{4}\phi^2+\frac{1}{4}\phi^3$. 
This, in turn, implies $(a,b)=(4,1)$ in section \ref{even P}, which does not provide with a 
solution of equation (\ref{eq4}). 

\vn
\underline{Case 2}: $\;(X,Y)=R$. We recall that $R=n_1Q_1+n_2Q_2$, with $n_1,n_2\in\Z$. 
In this case we are looking for points $(X,Y)=n_1Q_1+n_2Q_2$ with $X$ such that 
condition (\ref{X-condition 4.2}) be satisfied. More generally, we demand that  
the right-hand side of (\ref{inverse_bx+g}) be rational. For $(n_1,n_2)=(0,0)$ this condition 
is satisfied. Indeed, then $R={\mathcal O}$, 
$z(R)=0$ (by the definition of the function $z$; see section 4 of \cite{BT}), 
and the right-hand side of (\ref{inverse_bx+g}) is zero.

As mentioned immediately after (\ref{z_coord}), substitution of $z(R)$ in (\ref{inverse_bx+g}) 
from its value in (\ref{z_coord}) gives 
\[ \frac{1}{\beta X +\gamma} =
 \theta_0(n_1,n_2)+ \theta_1(n_1,n_2)\phi + \theta_2(n_1,n_2)\phi^2 + \theta_3(n_1,n_2)\phi^3 
\;,\]
hence, in order that the left-hand side be a rational number it is necessary that 
$\theta_1(n_1,n_2)=\theta_2(n_1,n_2)=\theta_3(n_1,n_2)=0$.
We will consider the system
\begin{equation} \label{system2}
\theta_3(n_1,n_2)=0,\: \theta_1(n_1,n_2)=0 \;\; n_1,n_2 \in\Z_3 \;,
\end{equation}
which, according to our discussion a few lines above, has the solution $(n_1,n_2)=(0,0)$,
and will show, using theorem \ref{a la Skolem} that this is its only solution in 3-adic 
integers. We set
\vspace{-3mm}
\[F_1(n_1,n_2)=\frac{1}{9}\theta_3(n_1,n_2)\,,\: F_2(n_1,n_2)=\frac{1}{9}\theta_1(n_1,n_2)\]
and we compute: 
\begin{align*}
F_1(n_1,n_2) =& 2n_1^2 +3(n_1^4+2n_1^2+n_1n_2) +3^2(\cdot) \\
F_2(n_1,n_2) = & n_1^2+n_1n_2+2n_2^2 +3(n_1^3n_2+n_1n_2^3+n_2^2+n_2^4) +3^2(\cdot) \;,
\end{align*}   
where $(\cdot)$ denotes an element of $\Z\langle n_1,n_2\rangle$ all of whose terms are of 
degree at least 2. Now, in the notation of theorem \ref{a la Skolem}, 
$f_{01}=2n_1^2, f_{02}=n_1^2+n_1n_2+2n_2^2$. We can obviously take 
$h_{11}=1, h_{21}=0, H_{1}=2n_1^2$. 
As for $H_2$, we can take it as the resultant of $f_{01},f_{02}$ with respect to $n_1$, 
finding thus $H_2=16n_2^4$  (here, $h_{12}=2n_2n_1-2n_2^2, h_{22}=-4n_2n_1+8n_2^2$, but we do not actually need these polynomials).
In view of the shape of the polynomials $H_1,H_2$, it follows by theorem \ref{a la Skolem} that 
$(n_1,n_2)=(0,0)$ is the only solution of $F_1(n_1,n_2)=0,F_2(n_1,n_2)=0$ in 3-adic integers 
and this solution corresponds to the zero point on the curve ${\mathcal E}_{10}$ which is 
of no interest for our initial problem. 

Summing up the previous results, we have proved the following
\begin{Results_r2} \label{Results_r2}
In the notation of section \ref{curve4}, the only points $(X,Y)$ on 
${\mathcal E}_{10}(\Q(\phi))$ satisfying the condition 
$\beta X+\gamma \in\Q$ ($\beta=6\phi+\phi^3, 
             \gamma=-4\phi-\phi^3$) are $\pm 2P_2, \pm (2P_1+2P_2)$.
No one of them furnishes a solution to equation (\ref{eq4}), hence no solution to our initial 
problem can be obtained from the elliptic curve ${\mathcal E}_{10}$.  
\end{Results_r2}

\newpage

{\small
\section {Appendix: The Mordell-Weil bases} \label{computations}
Notation: let $\nu$ be a non-Archimedean absolute value on $K$, where $K$ denotes 
$K_1$ or $K_2$, as appropriate, and let 
\[ \mbox{ord}_\nu: K_{\nu}^* \rightarrow {\bf Z} \]
be the corresponding normalized valuation: so that if the residue field
at $\nu$ has order $q_\nu$, then
\[ \log|x|_\nu = -\frac{1}{[K_\nu:{\bf Q}_\nu]} \mbox{ord}_\nu(x) \log(q_\nu) \]
for all $x \in K_\nu^*$. Equivalently,
\[ |x|_\nu^{[K_\nu:{\bf Q}_\nu]} = q_\nu^{-\mbox{ord}_\nu(x)}, \]
guaranteeing the product identity (over all non-Archimedean
and Archimedean absolute values)
\[ \prod_\nu |x|_\nu^{[K_\nu:{\bf Q}_\nu]} = 1. \]
The Archimedean valuation of $\Q$ has three extensions to $K$,
with $K_{\infty_1}=K_{\infty_2}={\bf R}$ and $K_{\infty_3}={\bf C}$.
We have $|x|_{\infty_1} = |x(\th)|$ (resp, $|x(\phi)|$), $|x|_{\infty_2} = |x(-\th)|$,
(resp. $|x(-\phi)|$), and $|x|_{\infty_3}=|x(i/\th)|$  (resp. $|x(2i/\ph)|$) - equivalently,
$|x|_{\infty_3}^2 = |x(i/\th)x(-i/\th)|$ (resp. $|x(2i/\ph)x(-2i/\ph)|$).\\
Define the indices $n_\nu = |K_\nu : \Q_\nu|$. Then
\[ \framebox{ $n_{(1+\th)}=4, \qquad n_\pi=4, \qquad n_{\infty_1}=n_{\infty_2}=1, \quad n_{\infty_3}=2$. } \]
The discriminants and Kodaira reduction types above $2$ 
are given in the following table; we also include the coefficients 
$\mu_{(1+\th)}$ and $\mu_\pi$, in Siksek's notation:\\ \\
\begin{tabular}{c|c|c|c|c|c}
Curve & Discriminant & Kodaira reduction type above $2$ & $\mu_{(1+\th)}$ & $\mu_\pi$  \\ \hline
(\ref{ell1}):  & $-\ep_1^{-14} \ep_2^6 (1+\th)^{18}$    &  $II$ & 0 & & \\
(\ref{ell2}):  & $-\ep_1^{-14} \ep_2^{12} (1+\th)^{18}$ &  $II$ & 0 & & \\
(\ref{ell3}):  & $-\ep_1^{-4} \ep_2^6 (1+\th)^{18}$     &  $II$ & 0 & & \\ 
(\ref{ell4}):  & $-\ep_1^{-4} \ep_2^{12} (1+\th)^{18}$  &  $II$ & 0 & & \\
(\ref{ell5}):  & $-\ep_1^{-2} \ep_2^{-12} \pi^{24}$     &  $I_4^*$ & & $1/4$ \\
(\ref{ell6}):  & $-\ep_1^4 \ep_2^{-12} \pi^{24}$        &  $I_6^*$ & & $1/4$ \\
(\ref{ell7}):  & $-\ep_1^{-2} \ep_2^{-6} \pi^{24}$      &  $I_4^*$ & & $1/4$ \\
(\ref{ell8}):  & $-\ep_1^4 \ep_2^{-6} \pi^{24}$         &  $I_6^*$ & & $1/4$ \\
(\ref{ell9}):  & $-\ep_1^2 \ep_2^{-12} \pi^{24}$        &  $I_4^*$ & & $1/4$ \\
(\ref{ell10}): & $-\ep_1^8 \ep_2^{-12} \pi^{24}$        &  $I_6^*$ & & $1/4$ \\
(\ref{ell11}): & $-\ep_1^2 \ep_2^{-6} \pi^{24}$         &  $I_4^*$ & & $1/4$ \\
(\ref{ell12}): & $-\ep_1^8 \ep_2^{-6} \pi^{24}$         &  $I_6^*$ & & $1/4$ \\ \hline
\end{tabular} \\ \\
We now make some remarks about the minimal polynomial of
$x(Q)$ for $Q \in E(K)$, with height $H(Q)$ bounded above by $B$, say.
Put $x_1 = x(Q)$.  If $|\Q(x_1) : \Q| = 4$, let $x_i$, $i=1,...,4$ 
denote the four conjugates of $x_1$, with minimum polynomial of 
$x_1$ being
\[ x^4 + a_1 x^3 + a_2 x^2 + a_3 x + a_4 = (x-x_1)(x-x_2)(x-x_3)(x-x_4). \]
Since
\[ |x_1|=|x_1|_{\infty_1} \leq \mbox{max}\{1,|x_1|_{\infty_1}\} \leq H(x_1)<B,\]
then
\[ |a_1|=|x_1+x_2+x_3+x_4| \leq |x_1|+|x_2|+|x_3|+|x_4| < 4 B, \]
using the fact that conjugate points have equal height.
In this way, we have
\begin{equation}
\label{coeffbounds4}
|a_1| < 4 B, \quad |a_2| < 6 B^2, \quad |a_3| < 4 B^3, \quad |a_4| < B^4.
\end{equation}
Similarly, if $|\Q(x_1) : \Q| = 2$, then the minimal polynomial
of $x_1$ is of type $x^2 + a_1 x + a_2$, where
\begin{equation}
\label{coeffbounds2}
|a_1| < 2 B, \qquad |a_2| < B^2.
\end{equation}
Finally, if $|\Q(x_1) : \Q| = 1$, then the minimal polynomial of
$x_1$ is of type $x + a_1$, where
\begin{equation}
\label{coeffbounds1}
|a_1| < B.
\end{equation}

\subsection{The curve $\mathcal{E}_1$ at (\ref{ell1})} \label{E1-computations}

\vn
>From the table of Kodaira reduction types, we have that (in Siksek's notation)
$\mu_\nu = 0$ except for
\[ \framebox{ $\mu_{\infty_1}=\mu_{\infty_2}=\mu_{\infty_3}=\frac{1}{3} .$ } \]
Further,
\[ \epsilon_\nu^{-1} = \inf_{(X,Y) \in E(K_\nu)} \frac{ \max(|f(X)|_\nu,|g(X)|_\nu)}{\max(1,|X|_\nu)^4} \]
with
\[ f(X)=  4 X^3 -4(\th+\th^2) X^2 + 4(1+\th+\th^3) X, \qquad g(X) = (X^2 -(1+\th+\th^3))^2. \]
Siksek gives a method for computing the $\epsilon_\nu$.
At $\infty_1$,
\[ \epsilon_{\infty_1}^{-1} = \inf_{(X,Y) \in E({\bf R}) } \frac{ \max(|f(X)|,|g(X)|)}{\max(1,|X|)^4}, \]
and the infimum needs to be taken over $X \in {\bf R}$ such that
$f(X) \geq 0$, that is, over
$[0, \infty)$. This infimum occurs at the root $4.275236449758861...$
of $f(X)-g(X)=0$, and has value $0.80190401917789682199...$,
so that
\[ \framebox{ $\epsilon_{\infty_1} = 1.24703203386508649515...$ }. \]
At $\infty_2$,
\[ \epsilon_{\infty_2}^{-1} = \inf_{(X,Y) \in E({\bf R}) } \frac{ \max(|\bar{f}(X)|,|\bar{g}(X)|)}{\max(1,|X|)^4}, \]
where
\[ \bar{f}(X) = 4 X^3 -4(-\th+\th^2) X^2 + 4(1-\th-\th^3) X,
 \quad \bar{g}(X) = (X^2 -(1-\th-\th^3))^2, \]
with infimum taken over $X \in {\bf R}$ such that
$f(X) \geq 0$, that is, over
$[0, \infty)$.
This infimum occurs at the root $0.021005066751861...$ of
$\bar{f}(X)-\bar{g}(X)=0$, and has value $0.00798861744730799360...$
so that 
\[ \framebox{ $\epsilon_{\infty_2} = 125.17810579814161228611...$ }. \]
At $\infty_3$, 
\[ \epsilon_{\infty_3}^{-1} = \inf_{(X,Y) \in E({\bf C})} \frac{\max(|F(X)|, |G(X)|)}{\max(1,|X|)^4} , \]
where
\[ F(X)^2 = 16 X^2 (X^4 +(4+2\th^2) X^3 +(9+3\th^2) X^2 +(10+4\th^2) X +(5+2\th^2)), \]
\[ G(X) = (X^4-2 X^2 +(5+2\th^2)). \] 
The infimum occurs at the root $ -1.45508613805...-0.5449200796689308...i$ of 
$|F(z)|=|G(z)|$, with value $0.6795900650263445248377698...$
(on the unit circle, the minimum taken exceeds $4$).
Thus
\[ \framebox{ $\epsilon_{\infty_3} = 1.471475307634514466025717... $ } \]
Putting the above together results in
\begin{eqnarray*}
h(P)-2 \hat{h}(P) & \leq & \frac{1}{4} (\frac{1}{3} \cdot 1 \cdot \log(1.24703203386508649515)\\
& & +\frac{1}{3} \cdot 1 \cdot \log(125.17810579814161228611)\\
& & +\frac{1}{3} \cdot 2 \cdot \log(1.4714753076345144660257)),
\end{eqnarray*}
that is,
\[ h(P)-2 \hat{h}(P) \leq 0.485252911746822... \]
Suppose now the point $G_1$ at (\ref{gen1}) is not a generator. 
We easily check that $G_1$ is not divisible by $2$ in $E(K)$, 
and so $G_1= m Q$ for $m \geq 3$ and $Q \in E(K)$. 
Note that since $(1+\th)^2 x(G_1) \in \mathcal{O}_K$, it follows that
$(1+\th)^2 x(Q) \in \mathcal{O}_K$.
Then
\[ h(Q) \leq 0.485252911746822+2 \hat{h}(Q) < 0.485252911746822+2 \hat{h}(G_1)/m^2 < 0.614184 \]
so that
\[ H(Q) < 1.84815. \]
Suppose first that $x(Q) \in \mathcal{O_K}$. Write $H(Q)<B$.
If $|\Q(x(Q)):\Q|=4$, then by direct computation, the minimum 
polynomial of $x(Q)$ is of type $X^4+4 c_1 X^3+2 c_2 X^2+4 c_3 X+c_4$, 
where $c_i \in \Z$, $i=1,..,4$. Similarly, if $|\Q(x(Q)):\Q|=2$,
then the minimal polynomial of $x(Q)$ is of type $X^2+2 c_1 X+c_2$,
with $c_i \in \Z$. From (\ref{coeffbounds4}), (\ref{coeffbounds2}),
(\ref{coeffbounds1}), we therefore have to investigate the following 
polynomials:\\
$\bullet$ $X^4+4 c_1 X^3+2 c_2 X^2 +4 c_3 X+ c_4$, $c_i \in \Z$,
$|c_1| < B$, $|c_2| < 3 B^2$, $|c_3| < B^3$, $|c_4| < B^4$.\\
$\bullet$ $X^2+2 c_1 X + c_2$,  $c_i \in \Z$, $|c_1| < B$, $|c_2| < B^2$. \\
$\bullet$ $X+c_1$, $c_1 \in \Z$, $|c_1| < B$.\\ \\
Suppose second that $x(Q) = u/(1+\th)^2$, where $u \in \mathcal{O}_K$,
and $u \equiv 1 \pmod {(1+\th)}$. If $|\Q(x(Q)):\Q|=4$, then by direct 
computation, the minimum polynomial of $x(Q)$ is of type 
$X^4+4 c_1 X^3+c_2 X^2+2 c_3 X+\frac{c_4}{4}$, where $c_i \in \Z$, 
and $c_2 \equiv c_4 \equiv 1 \pmod 2$. Similarly, if
$|\Q(x(Q)):\Q|=2$, then the minimal polynomial of $x(Q)$ is
of type $X^2 + 2 c_1 X+\frac{c_2}{4}$, where $c_i \in \Z$, and
$c_2 \equiv 3 \pmod 4$.
As above, we then have to investigate polynomials:\\
$\bullet$ $X^4+4 c_1 X^3+c_2 X^2+2 c_3 X+\frac{c_4}{4}$, $c_i \in \Z$, $c_2 \equiv c_4 \equiv 1 \pmod 2$, 
\[ |c_1| < B, |c_2| < 6 B^2, |c_3| < 2 B^3, |c_4| < 4 B^4. \]
$\bullet$ $X^2 + 2c_1 X+\frac{c_2}{4}$, $c_i \in \Z$, $c_2 \equiv 3 \pmod 4$, $|c_1| < B$, $|c_2| < 4 B^2$.\\ \\
Numerically, we have to investigate polynomials:\\
$\bullet$ $X^4+4 c_1 X^3+2 c_2 X^2 +4 c_3 X+ c_4$, $c_i \in \Z$,
$|c_1| \leq 1$, $|c_2| \leq 10$, $|c_3| \leq 6$, $|c_4| \leq 11$ \\
$\bullet$ $X^2+2 c_1 X + c_2$,  $c_i \in \Z$, $|c_1| \leq 1$, $|c_2| \leq 3$ \\
$\bullet$ $X+c_1$, $c_1 \in \Z$, $|c_1| \leq 1$\\
and\\
$\bullet$ $X^4+4 c_1 X^3+c_2 X^2+2 c_3 X+\frac{c_4}{4}$, $c_i \in \Z$, $c_2 \equiv c_4 \equiv 1 \pmod 2$,
$|c_1| \leq 1$, $|c_2| \leq 20$, $|c_3| \leq 12$, $|c_4| \leq 46$ \\
$\bullet$ $X^2 + 2c_1 X+\frac{c_2}{4}$, $c_i \in \Z$, $c_2 \equiv 3 \pmod 4$, $|c_1| \leq 1$, $|c_2| < 13$.\\
Each polynomial has to be tested to see if a root can be the $X$-coordinate
of a point in $\mathcal{E}_1(K)$.
Computation shows that in the given range, only the points  $\pm G_1$
arise.
It follows that $G_1$ is indeed
a generator of the group of points defined over $K$.

\subsection{The curve $\mathcal{E}_2$ at (\ref{ell2})} \label{E2-computations}

\vn
>From the table of Kodaira reduction types, we have $\mu_\nu = 0$ except for
\[ \framebox{ $\mu_{\infty_1}=\mu_{\infty_2}=\mu_{\infty_3}=\frac{1}{3} .$ } \]
Further,
\[ \epsilon_\nu^{-1} = \inf_{(X,Y) \in E(K_\nu)} \frac{ \max(|f(X)|_\nu,|g(X)|_\nu)}{\max(1,|X|_\nu)^4} \]
with
\[ f(X)= 4 X^3 -4(\th-\th^2) X^2 + 4(1-\th-\th^3) X, \quad g(X) = (X^2 -(1-\th-\th^3))^2. \]
At $\infty_1$,
\[ \epsilon_{\infty_1}^{-1} = \inf_{(X,Y) \in E({\bf R}) } \frac{ \max(|f(X)|,|g(X)|)}{\max(1,|X|)^4}, \]
with infimum taken over $X \in \R$ such that
$f(X) \geq 0$, that is, over
$[0, \infty)$. This infimum occurs at the root $0.023441018652769...$
of $f(X)-g(X)=0$, and has value $0.00796927528986859148...$,
so that
\[ \framebox{ $\epsilon_{\infty_1} = 125.48192446950711297112...$ }. \]
At $\infty_2$,
\[ \epsilon_{\infty_2}^{-1} = \inf_{(X,Y) \in E({\bf R}) } \frac{ \max(|\bar{f}(X)|,|\bar{g}(X)|)}{\max(1,|X|)^4}, \]
where
\[ \bar{f}(X) = 4 X^3 +4(\th+\th^2) X^2 + 4(1+\th+\th^3) X, \quad \bar{g}(X) = (X^2 -(1+\th+\th^3))^2, \]
with infimum taken over $X \in {\bf R}$ such that $f(X) \geq 0$, that is, 
over $[0, \infty)$.
This infimum occurs at the root $5.645614058038130...$ of
$\bar{f}(X)-\bar{g}(X)=0$, and has value $0.88372963806597132831...$
so that
\[ \framebox{ $\epsilon_{\infty_2} = 1.13156779735087822111...$ }. \]
At $\infty_3$,
\[ \epsilon_{\infty_3}^{-1} = \inf_{(X,Y) \in E({\bf C})} \frac{\max(|F(X)|, |G(X)|)}{\max(1,|X|)^4}, \]
where
\[ F(X)^2 = 16 X^2 (X^4 -(4+2\th^2)X^3 + (9+3\th^2)X^2 -(10+4\th^2)X + (5+2\th^2)), \]
\[ G(X) =  (X^4-2 X^2 +(5+2\th^2)). \]
The infimum occurs at the root $1.455086138050493956497..-0.544920079668930802096...i$ of
$|F(z)|=|G(z)|$, with value $0.6795900650263445248377698...$
(on the unit circle, the minimum taken exceeds $4$).
Thus
\[ \framebox{ $\epsilon_{\infty_3} = 1.471475307634514466025717... $ } \]
Putting the above together results in
\begin{eqnarray*}
h(P)-2 \hat{h}(P) & \leq & \frac{1}{4} (\frac{1}{3} \cdot 1 \cdot \log(125.48192446950711297112)\\
 & & + \frac{1}{3} \cdot 1 \cdot \log(1.13156779735087822111) \\
 & & + \frac{1}{3} \cdot 2 \cdot \log(1.471475307634514466025717)),
\end{eqnarray*}
that is,
\[ h(P)-2 \hat{h}(P) \leq   0.477358069897830... \]
Suppose now the point $G_2$ at (\ref{gen2}) is not a generator.
We easily check that $G_2$ is not divisible by $2$ in $E(K)$,
and so $G_2= m Q$ for $m \geq 3$ and $Q \in E(K)$.
Note that since $(1+\th)^2 x(G_2) \in \mathcal{O}_K$, it follows that
$(1+\th)^2 x(Q) \in \mathcal{O}_K$.
Then
\[ h(Q) \leq 0.47735806989783+2 \hat{h}(Q) < 0.47735806989783+2*\hat{h}(G_2)/m^2 < 0.533699 \]
so that
\[ H(Q) < 1.70523. \]
Arguing as in the previous instance, we have to consider all
polynomials of type\\
$\bullet$ $X^4+4 c_1 X^3+2 c_2 X^2 +4 c_3 X+ c_4$, $c_i \in \Z$,
$|c_1| \leq 1$, $|c_2| \leq 8$, $|c_3| \leq 4$, $|c_4| \leq 8$ \\
$\bullet$ $X^2+2 c_1 X + c_2$,  $c_i \in \Z$, $|c_1| \leq 1$, $|c_2| \leq 2$ \\
$\bullet$ $X+c_1$, $c_1 \in \Z$, $|c_1| \leq 1$\\
and\\
$\bullet$ $X^4+4 c_1 X^3+c_2 X^2+2 c_3 X+\frac{c_4}{4}$, $c_i \in \Z$, $c_2 \equiv c_4 \equiv 1 \pmod 2$,
$|c_1| \leq 1$, $|c_2| \leq 17$, $|c_3| \leq 9$, $|c_4| \leq 33$ \\
$\bullet$ $X^2 + 2c_1 X+\frac{c_2}{4}$, $c_i \in \Z$, $c_2 \equiv 3 \pmod 4$, $|
c_1| \leq 1$, $|c_2| < 11$.\\
Each polynomial has to be tested to see if a root can be the $X$-coordinate
of a point in $\mathcal{E}_2(K)$.
Computation shows that in the given range, only the points $\pm G_2$, 
$\pm G_2 +(0,0)$ arise.
It follows that $G_2$ is indeed
a generator of the group of points defined over $K$.

\subsection{The curve $\mathcal{E}_3$ at (\ref{ell3})} \label{E3-computations}

>From the table of Kodaira reduction types, we have $\mu_\nu = 0$ except for
\[ \framebox{ $\mu_{\infty_1}=\mu_{\infty_2}=\mu_{\infty_3}=\frac{1}{3} .$ } \]
Further,
\[ \epsilon_\nu^{-1} = \inf_{(X,Y) \in E(K_\nu)} \frac{ \max(|f(X)|_\nu,|g(X)|_\nu)}{\max(1,|X|_\nu)^4} \]
with
\[ f(X) = 4 X^3 -4(1+\th) X^2 + 4(\th+\th^2-\th^3) X, \quad g(X) = (X^2 -(\th+\th^2-\th^3))^2. \]
At $\infty_1$,
\[ \epsilon_{\infty_1}^{-1} = \inf_{(X,Y) \in E({\bf R}) } \frac{ \max(|f(X)|,|g(X)|)}{\max(1,|X|)^4}, \]
and the infimum needs to be taken over $X \in {\bf R}$ such that
$f(X) \geq 0$, that is, over
$[0, \infty)$. This infimum occurs at the turning point $0.738691905746190...$
of $f(X)0$, and has value $0.36278136846310610700...$,
so that
\[ \framebox{ $\epsilon_{\infty_1} = 2.75648113969143186636...$ }. \]
At $\infty_2$,
\[ \epsilon_{\infty_2}^{-1} = \inf_{(X,Y) \in E({\bf R}) } \frac{ \max(|\bar{f}(X)|,|\bar{g}(X)|)}{\max(1,|X|)^4}, \]
where
\[ \bar{f}(X) = 4 X^3 -4(1-\th) X^2 + 4(-\th+\th^2+\th^3) X, \quad \bar{g}(X) = (X^2 -(-\th+\th^2+\th^3))^2, \]
with infimum taken over $X \in {\bf R}$ such that $f(X) \geq 0$, that is, over
$[0, \infty)$.
This infimum occurs at the root $0.010221121380833...$ of
$\bar{f}(X)-\bar{g}(X)=0$, and has value $0.00137643273231028235...$
so that
\[ \framebox{ $\epsilon_{\infty_2} = 726.51570725257570965658...$ }. \]
At $\infty_3$,
\[ \epsilon_{\infty_3}^{-1} = \inf_{(X,Y) \in E({\bf C})} \frac{\max(|F(X)|, |G(X)|)}{\max(1,|X|)^4}, \]
where
\[ F(X)^2 = 16 X^2 (X^4 - 2 X^3 -(1+\th^2) X^2 -(10+4\th^2) X + (29+12\th^2)), \]
\[ G(X) =  (X^4-2 X^2 +(5+2\th^2)). \]
The infimum occurs at the root $-1.37342963506048574888...-1.985476809807611687126...i$ of
$|F(z)|=|G(z)|$, with value $0.1229849339729954943136161149...$
(on the unit circle, the minimum taken exceeds $24$).
Thus
\[ \framebox{ $\epsilon_{\infty_3} = 8.13107726040309719656378512... $ } \]
Putting the above together results in
\begin{eqnarray*}
h(P)-2 \hat{h}(P) & \leq & \frac{1}{4} (\frac{1}{3} \cdot 1 \cdot \log(2.75648113969143186636)\\
 & & + \frac{1}{3} \cdot 1 \cdot \log(726.51570725257570965658) \\
 & & + \frac{1}{3} \cdot 2 \cdot \log(8.13107726040309719656378512)),
\end{eqnarray*}
that is,
\[ h(P)-2 \hat{h}(P) \leq  0.982800154866326... \]
Suppose now the point $G_3$ at (\ref{gen3}) is not a generator.
We easily check that $G_3$ is not divisible by $2$ in $E(K)$,
and so $G_3= m Q$ for $m \geq 3$ and $Q \in E(K)$.
Note that since $(1+\th)^2 x(G_3) \in \mathcal{O}_K$, it follows that
$(1+\th)^2 x(Q) \in \mathcal{O}_K$.
Then
\[ h(Q) \leq 0.982800154866326 + 2 \hat{h}(Q) < 0.982800154866326 + 2 \hat{h}(G_3)/m^2 < 1.037355 \]
so that
\[ H(Q) < 2.82175. \]
Arguing as in the previous instance, we have to consider all
polynomials of type\\
$\bullet$ $X^4+4 c_1 X^3+2 c_2 X^2 +4 c_3 X+ c_4$, $c_i \in \Z$,
$|c_1| \leq 2$, $|c_2| \leq 23$, $|c_3| \leq 22$, $|c_4| \leq 63$ \\
$\bullet$ $X^2+2 c_1 X + c_2$,  $c_i \in \Z$, $|c_1| \leq 2$, $|c_2| \leq 7$ \\
$\bullet$ $X+c_1$, $c_1 \in \Z$, $|c_1| \leq 2$\\
and\\
$\bullet$ $X^4+4 c_1 X^3+c_2 X^2+2 c_3 X+\frac{c_4}{4}$, $c_i \in \Z$, $c_2 \equiv c_4 \equiv 1 \pmod 2$,
$|c_1| \leq 2$, $|c_2| \leq 47$, $|c_3| \leq 44$, $|c_4| \leq 253$ \\
$\bullet$ $X^2 + 2c_1 X+\frac{c_2}{4}$, $c_i \in \Z$, $c_2 \equiv 3 \pmod 4$, $|
c_1| \leq 2$, $|c_2| < 31$.\\
Each polynomial has to be tested to see if a root can be the $X$-coordinate
of a point in $\mathcal{E}_3(K)$.
Computation shows that in the given range, only the points $\pm G_3$,
$\pm G_3+(0,0)$ arise.
It follows that $G_3$ is indeed
a generator of the group of points defined over $K$.

\subsection{The curve $\mathcal{E}_4$ at (\ref{ell4})} \label{E4-computations}

>From the table of Kodaira reduction types, we have $\mu_\nu = 0$ except for
\[ \framebox{ $\mu_{\infty_1}=\mu_{\infty_2}=\mu_{\infty_3}=\frac{1}{3} .$ } \]
Further,
\[ \epsilon_\nu^{-1} = \inf_{(X,Y) \in E(K_\nu)} \frac{ \max(|f(X)|_\nu,|g(X)|_\nu)}{\max(1,|X|_\nu)^4} \]
with
\[ f(X) = 4 X^3 -4(1-\th) X^2 + 4(-\th+\th^2+\th^3) X, \quad g(X) = (X^2 -(-\th+\th^2+\th^3))^2. \]
The curve is the conjugate of the curve (\ref{ell3}) under
$\th \rightarrow -\th$, and so 
\[ \framebox{ $\epsilon_{\infty_1} = 726.51570725257570965658...$ }, \]
\[ \framebox{ $\epsilon_{\infty_2} = 2.75648113969143186636...$ } \]
At $\infty_3$,
\[ \epsilon_{\infty_3}^{-1} = \inf_{(X,Y) \in E({\bf C})} \frac{\max(|F(X)|, |G(X)|)}{\max(1,|X|)^4}, \]
where
\[ F(X)^2 =  16 X^2 (X^4 - 2 X^3 -(1+\th^2) X^2 -(10+4\th^2) X + (29+12\th^2)), \]
\[ G(X) =  (X^4-2 X^2 +(5+2\th^2)). \]
The infimum occurs at the root $-1.37342963506048574888049...-1.985476809807611687126...i$
of $|F(z)|=|G(z)|$, with value $0.1229849339729954943136161149...$
(on the unit circle, the minimum taken exceeds $24$).
Thus
\[ \framebox{ $\epsilon_{\infty_3} = 8.13107726040309719656378512...$ } \]
Putting the above together results in
\begin{eqnarray*}
h(P)-2 \hat{h}(P) & \leq & \frac{1}{4} (\frac{1}{3} \cdot 1 \cdot \log(726.51570725257570965658)\\
 & & + \frac{1}{3} \cdot 1 \cdot \log(2.75648113969143186636) \\
 & & + \frac{1}{3} \cdot 2 \cdot \log(8.13107726040309719656378512)),
\end{eqnarray*}
that is,
\[ h(P)-2 \hat{h}(P) \leq 0.982800154866326... \]
Suppose now the point $G_4$ at (\ref{gen4}) is not a generator.
We easily check that $G_4$ is not divisible by $2$ in $E(K)$,
and so $G_4= m Q$ for $m \geq 3$ and $Q \in E(K)$.
Note that since $(1+\th)^2 x(G_4) \in \mathcal{O}_K$, it follows that
$(1+\th)^2 x(Q) \in \mathcal{O}_K$.
Then
\[ h(Q) \leq 0.982800154866326 + 2 \hat{h}(Q) < 0.982800154866326 + 2 \hat{h}(G_4)/m^2 < 1.037355 \]
so that
\[ H(Q) < 2.82175. \]
Arguing as in the previous instance, we have to consider all
polynomials of type\\
$\bullet$ $X^4+4 c_1 X^3+2 c_2 X^2 +4 c_3 X+ c_4$, $c_i \in \Z$,
$|c_1| \leq 2$, $|c_2| \leq 23$, $|c_3| \leq 22$, $|c_4| \leq 63$ \\
$\bullet$ $X^2+2 c_1 X + c_2$,  $c_i \in \Z$, $|c_1| \leq 2$, $|c_2| \leq 7$ \\
$\bullet$ $X+c_1$, $c_1 \in \Z$, $|c_1| \leq 2$\\
and\\
$\bullet$ $X^4+4 c_1 X^3+c_2 X^2+2 c_3 X+\frac{c_4}{4}$, $c_i \in \Z$, $c_2 \equiv c_4 \equiv 1 \pmod 2$,
$|c_1| \leq 2$, $|c_2| \leq 47$, $|c_3| \leq 44$, $|c_4| \leq 253$ \\
$\bullet$ $X^2 + 2c_1 X+\frac{c_2}{4}$, $c_i \in \Z$, $c_2 \equiv 3 \pmod 4$, $|
c_1| \leq 2$, $|c_2| < 31$.\\
Each polynomial has to be tested to see if a root can be the $X$-coordinate
of a point in $\mathcal{E}_4(K)$.
Computation shows that in the given range, only the points $G_4$,
$\pm G_4+(0,0)$ arise.
It follows that $G_4$ is indeed
a generator of the group of points defined over $K$.

\subsection{The curve $\mathcal{E}_5$ at (\ref{ell5})} \label{E5-computations}

>From the table of Kodaira reductions, we have $\mu_\nu = 0$ except for
\[ \framebox{ $\mu_{\pi}=\frac{1}{4}, \qquad \mu_{\infty_1}=\mu_{\infty_2}=\mu_{\infty_3}=\frac{1}{3}.$ } \]
Further,
\[ \epsilon_\nu^{-1} = \inf_{(X,Y) \in E(K_\nu)} \frac{ \max(|f(X)|_\nu,|g(X)|_\nu)}{\max(1,|X|_\nu)^4} \]
with
\[ f(X) = 4 X^3 -4\ph X^2 + (4+2\ph^2) X, \quad g(X) = (X^2 -(1+\ph^2/2))^2. \]
Siksek gives a method for computing the $\epsilon_\nu$. For the
non-Archimedean valuation, we have the following (in Siksek's notation).
First, we observe that $g(1-\frac{1}{2}\phi+\frac{1}{4}\phi^3) \equiv 0 \pmod {\pi^{10}}$,
and $g(X) \not \equiv 0 \pmod {\pi^{12}}$ for any $X \in K$. Thus 
$\ep_\pi = |\pi|_\pi^{-2 j} = (2^{-\frac{1}{4}})^{-2 j}$, where $j \leq 5$. 
This weak inequality is all that we need, resulting in
\[ \framebox{ $\epsilon_{\pi} \leq 2^{\frac{5}{2}}$}. \]
At $\infty_1$,
\[ \epsilon_{\infty_1}^{-1} = \inf_{(X,Y) \in E({\bf R}) } \frac{ \max(|f(X)|,|g(X)|)}{\max(1,|X|)^4}, \]
and the infimum needs to be taken over $X \in {\bf R}$ such that
$f(X) \geq 0$, that is, over
$[0, \infty)$. This infimum occurs at the $4.108570541436509...$
of $f(X)=g(X)$, and has value $0.83946151126494434491...$,
so that
\[ \framebox{ $\epsilon_{\infty_1} = 1.19123984432966783131...$ }. \]
At $\infty_2$,
\[ \epsilon_{\infty_2}^{-1} = \inf_{(X,Y) \in E({\bf R}) } \frac{ \max(|\bar{f}(X)|,|\bar{g}(X)|)}{\max(1,|X|)^4}, \]
where
\[ \bar{f}(X) = 4 X^3 + 4\ph X^2 + (4+2\ph^2) X, \quad \bar{g}(X) = (X^2 -(1+\ph^2/2))^2, \]
with infimum taken over $X \in {\bf R}$ such that
$f(X) \geq 0$, that is, over
$[0, \infty)$.
This infimum occurs at the root $5.383909674320621...$ of
$\bar{f}(X)-\bar{g}(X)=0$, and has value $0.90480288995171512682...$
so that
\[ \framebox{ $\epsilon_{\infty_2} = 1.10521309238232547422...$ }. \]
At $\infty_3$,
\[ \epsilon_{\infty_3}^{-1} = \inf_{(X,Y) \in E({\bf C})} \frac{\max(|F(X)|, |G(X)|)}{\max(1,|X|)^4}, \]
where 
\[ F(X)^2 = 16 X^2 (X^4 + 2 X^2 +2), \]
\[ G(X) = X^4  + (2+\ph^2) X^2 + 2. \]
The infimum occurs at the root $-0.444261439847776944198...-1.103107127815551338132621...i$ of
$|F(z)|=|G(z)|$, with value $0.5582416466277690341698809...$
(on the unit circle, the minimum taken exceeds $2$).
Thus
\[ \framebox{ $\epsilon_{\infty_3} = 1.791338940834688072056363...$ } \]
Putting the above together results in
\begin{eqnarray*}
h(P)-2 \hat{h}(P) & \leq & \frac{1}{4} (\frac{1}{4} \cdot 4 \cdot \log 2^{\frac{5}{2}} + \frac{1}{3} \cdot 1 \cdot \log(1.19123984432966783131)\\
 & & + \frac{1}{3} \cdot 1 \cdot \log(1.10521309238232547422)\\
 & & + \frac{1}{3} \cdot 2 \cdot \log(1.791338940834688072056363)),
\end{eqnarray*}
that is,
\[ h(P)-2 \hat{h}(P) \leq 0.553296947402687... \]
Suppose now the point $G_5$ at (\ref{gen5}) is not a generator.
We easily check that $G_5$ is not divisible by $2$ in $E(K)$,
and so $G_5= m Q$ for $m \geq 3$ and $Q \in E(K)$, with 
$x(Q) \in \mathcal{O}_K$.
Then
\[ h(Q) \leq 0.553296947402687 + 2 \hat{h}(Q) < 0.553296947402687 + 2 \hat{h}(G_5)/m^2 < 0.609176 \]
so that
\[ H(Q) < 1.83892. \]
Write $H(Q) < B$. By direct computation, if $|\Q(x(Q)):\Q|=4$, then 
the minimal polynomial 
for x(Q) is of type $X^4 +4 a_1 X^3 + 2 a_2 X^2 + 4 a_3 X + a_4$, 
with $a_i \in \Z$, and, from (\ref{coeffbounds4}), $|a_1| < B$, $|a_2| < 3 B^2$,
$|a_3| < B^3$, $|a_4| < B^4$. Similarly, if $|\Q(x(Q)):\Q|=2$, then the minimal
polynomial of $x(Q)$ is of type $X^2 +2 a_1 X +a_2$, $a_i \in Z$, with, from 
(\ref{coeffbounds2}), $|a_1| < B$, $|a_2| < B^2$.
Accordingly, 
we have to consider polynomials of the following types, where $a_i \in \Z$:
\begin{eqnarray*}
x^4 + 4 a_1 x^3 + 2 a_2 x^2 + 4 a_3 x + a_4, & |a_1| \leq 1, |a_2| \leq 10, |a_3| \leq 6, |a_4| \leq 11, & \\
x^2+2 a_1 x+a_2, & |a_1| \leq 1, |a_2| \leq 3, & \\
x+a_1, & |a_1| \leq 1. &
\end{eqnarray*}
Each polynomial has to be tested to see if a root can be the $X$-coordinate
of a point in $\mathcal{E}_5(K)$.
Computation shows that in the given range, only the points $\pm G_5$,
$\pm G_5+(0,0)$ arise.
It follows that $G_5$ is indeed
a generator of the group of points defined over $K$.

\subsection{The curve $\mathcal{E}_6$ at (\ref{ell6})} \label{E6-computations}

>From the table of Kodaira reductions, we have $\mu_\nu = 0$ except for
\[ \framebox{ $\mu_{\pi}=\frac{1}{4}, \qquad \mu_{\infty_1}=\mu_{\infty_2}=\mu_{\infty_3}=\frac{1}{3}.$ } \]
Further,
\[ \epsilon_\nu^{-1} = \inf_{(X,Y) \in E(K_\nu)} \frac{ \max(|f(X)|_\nu,|g(X)|_\nu)}{\max(1,|X|_\nu)^4} \]
with
\[ f(X) = 4 X^3 +(-4+2\ph^2) X^2 + (4-2\ph^2) X, \quad g(X) = (X^2 -(1-\ph^2/2))^2. \]
Siksek gives a method for computing the $\epsilon_\nu$. For the
non-Archimedean valuation, we have the following (in Siksek's notation).
First, we observe that 
$g(1-\frac{1}{2}\phi-\frac{1}{4}\phi^3) \equiv 0 \pmod {\pi^{10}}$,
and $g(X) \not \equiv 0 \pmod {\pi^{12}}$ for any $X \in K$. Thus
$\ep_\pi = |\pi|_\pi^{-2 j} = (2^{-\frac{1}{4}})^{-2 j}$, where $j \leq 5$.
This weak inequality is all that we need, resulting in
\[ \framebox{ $\epsilon_{\pi} \leq 2^{\frac{5}{2}}$}. \]
At $\infty_1$,
\[ \epsilon_{\infty_1}^{-1} = \inf_{(X,Y) \in E({\bf R}) } \frac{ \max(|f(X)|,|g(X)|)}{\max(1,|X|)^4}, \]
and the infimum needs to be taken over $X \in {\bf R}$ such that
$f(X) \geq 0$, that is, over
$[0, \infty)$. This infimum occurs at the root $0.152240934977426...$
of $f(X)=g(X)$, and has value $0.31652903917264027803...$,
so that
\[ \framebox{ $\epsilon_{\infty_1} = 3.15926779613602254445...$ }. \]
At $\infty_2$, the curve is invariant under $\ph \rightarrow -\ph$,
and so
\[ \framebox{ $\epsilon_{\infty_2} = 3.15926779613602254445...$ }. \]
At $\infty_3$,
\[ \epsilon_{\infty_3}^{-1} = \inf_{(X,Y) \in E({\bf C})} \frac{\max(|F(X)|, |G(X)|)}{\max(1,|X|)^4}, \]
where
\[ F(X)^2 = 16 X^2 (X^4 - (6+\ph^2) X^3 +(16+3\ph^2) X^2 -(20+4\ph^2) X + (10+2\ph^2)), \]
\[ G(X) =  X^4 - (6+\ph^2) X^2 + (10+2\ph^2). \] 
The infimum occurs at the root $1.797932651931813404063...-0.426206219441401112133512...i$ of
$|F(z)|=|G(z)|$, with value $0.212818253072924089328469775...$
(on the unit circle, the minimum taken exceeds $4$).
Thus
\[ \framebox{ $\epsilon_{\infty_3} = 4.6988450734878506972817404... $ } \]
Putting the above together results in
\begin{eqnarray*}
h(P)-2 \hat{h}(P) & \leq & \frac{1}{4} (\frac{1}{4} \cdot 4 \cdot \log 2^{\frac{5}{2}} + \frac{1}{3} \cdot 1 \cdot \log(3.15926779613602254445)\\
 & & + \frac{1}{3} \cdot 1 \cdot \log(3.15926779613602254445) \\ 
 & & + \frac{1}{3} \cdot 2 \cdot \log(4.69884507348785069)),
\end{eqnarray*}
that is,
\[ h(P)-2 \hat{h}(P) \leq 0.882826494540115... \]
Suppose now the point $G_6$ at (\ref{gen6}) is not a generator.
We easily check that $G_6$ is not divisible by $2$ in $E(K)$,
and so $G_6=m Q$ for $m \geq 3$ and $Q \in E(K)$, with
$x(Q) \in \mathcal{O}_K$.
Then
\[ h(Q) \leq 0.882826494540115 + 2 \hat{h}(Q) < 0.882826494540115 + 2 \hat{h}(G_6)/m^2 < 0.923750 \]
so that
\[ H(Q) < 2.51872. \]
Arguing as in the case of the curve (\ref{ell5}), we must
consider polynomials of the following types, where $a_i \in \Z$:
\begin{eqnarray*}
x^4 + 4 a_1 x^3 + 2 a_2 x^2 + 4 a_3 x + a_4, & |a_1| \leq 2, |a_2| \leq 19, |a_3| \leq 15, |a_4| \leq 40, & \\
x^2+2 a_1 x+a_2, & |a_1| \leq 2, |a_2| \leq 6, & \\
x+a_1, & |a_1| \leq 2. &
\end{eqnarray*}
Each polynomial has to be tested to see if a root can be the $X$-coordinate
of a point in $\mathcal{E}_6(K)$.
Computation shows that in the given range, only the points $\pm G_6$,
$\pm G_6+(0,0)$ arise.
It follows that $G_6$ is indeed
a generator of the group of points defined over $K$.


\subsection{The curve $\mathcal{E}_7$ at (\ref{ell7})} \label{E7-computations}
>From the table of Kodaira reductions, we have $\mu_\nu = 0$ except for
\[ \framebox{ $\mu_{\pi}=\frac{1}{4}, \qquad \mu_{\infty_1}=\mu_{\infty_2}=\mu_{\infty_3}=\frac{1}{3}.$ } \]
Further,
\[ \epsilon_\nu^{-1} = \inf_{(X,Y) \in E(K_\nu)} \frac{ \max(|f(X)|_\nu,|g(X)|_\nu)}{\max(1,|X|_\nu)^4} \]
with
\begin{eqnarray*} 
f(X) & = & 4 X^3 +(-8-8\ph-2 \ph^3) X^2 + (52+56\ph+10\ph^2+12\ph^3) X, \\
g(X) & = & (X^2 -(13+14\ph+\frac{5}{2}\ph^2+3\ph^3))^2.
\end{eqnarray*}
Siksek gives a method for computing the $\epsilon_\nu$. For the
non-Archimedean valuation, we have the following (in Siksek's notation).
First, we observe that 
$g(1+\frac{1}{2}\phi-\frac{1}{4}\phi^3) \equiv 0 \pmod {\pi^{10}}$,
and $g(X) \not \equiv 0 \pmod {\pi^{12}}$ for any $X \in K$. Thus
$\ep_\pi = |\pi|_\pi^{-2 j} = (2^{-\frac{1}{4}})^{-2 j}$, where $j \leq 5$.
This weak inequality is all that we need, resulting in
\[ \framebox{ $\epsilon_{\pi} \leq 2^{\frac{5}{2}}$}. \]
At $\infty_1$,
\[ \epsilon_{\infty_1}^{-1} = \inf_{(X,Y) \in E({\bf R}) } \frac{ \max(|f(X)|,|g(X)|)}{\max(1,|X|)^4}, \]
and the infimum needs to be taken over $X \in {\bf R}$ such that
$f(X) \geq 0$, that is, over
$[0, \infty)$. This infimum occurs at the root $9.043006133337668...$
of $f(X)=g(X)$, and has value $ 0.39970098305719519573...$,
so that
\[ \framebox{ $\epsilon_{\infty_1} = 2.50187025398660338324...$ }. \]
At $\infty_2$,
\[ \epsilon_{\infty_2}^{-1} = \inf_{(X,Y) \in E({\bf R}) } \frac{ \max(|\bar{f}(X)|,|\bar{g}(X)|)}{\max(1,|X|)^4}, \]
where
\begin{eqnarray*}
\bar{f}(X) & = & 4 X^3 + (-8+8\ph+2\ph^3) \ph X^2 + (52-56\ph+10\ph^2-12\ph^3) X, \\
\bar{g}(X) & = & (X^2 -(13-14\ph+\frac{5}{2}\ph^2-3\ph^3))^2,
\end{eqnarray*}
with infimum taken over $X \in {\bf R}$ such that
$f(X) \geq 0$, that is, over
$[0, \infty)$.
This infimum occurs at the root $0.015710679827598...$ of
$\bar{f}(X)-\bar{g}(X)=0$, and has value $0.00438935169160511858...$
so that
\[ \framebox{ $\epsilon_{\infty_2} = 227.82407750842934587031...$ }. \]
At $\infty_3$,
\[ \epsilon_{\infty_3}^{-1} = \inf_{(X,Y) \in E({\bf C})} \frac{\max(|F(X)|, |G(X)|)}{\max(1,|X|)^4}, \]
where
\[ F(X)^2 = 16 X^2 (X^4 -4 X^3 +(10-4\ph^2) X^2 -(4+2\ph^2) X + 2), \]
\[ G(X) =  X^4 +(-6+5\ph^2) X^2 + 2. \]
The infimum occurs at the root $1.164435178539534874799...-0.24146278668295160003697...i$ of
$|F(z)|=|G(z)|$, with value $0.52138210146214758954528399...$
(on the unit circle, the minimum taken exceeds $0.90$).
Thus
\[ \framebox{ $\epsilon_{\infty_3} = 1.9179791504074102437227773... $ } \]
Putting the above together results in
\begin{eqnarray*}
h(P)-2 \hat{h}(P) & \leq & \frac{1}{4} (\frac{1}{4} \cdot 4 \cdot \log 2^{\frac{5}{2}} + \frac{1}{3} \cdot 1 \cdot \log(2.50187025398660338324)\\
 & & + \frac{1}{3} \cdot 1 \cdot \log(227.82407750842934587031) \\
 & & + \frac{1}{3} \cdot 2 \cdot \log(1.917979150407410243722)),
\end{eqnarray*}
that is,
\[ h(P)-2 \hat{h}(P) \leq 1.070563363421848... \]
Suppose now the point $G_7$ at (\ref{gen7}) is not a generator.
We easily check that $G_7$ is not divisible by $2$ in $E(K)$,
and so $G_7=m Q$ for $m \geq 3$ and $Q \in E(K)$, with
$x(Q) \in \mathcal{O}_K$.
Then
\[ h(Q) \leq 1.070563363421848 + 2 \hat{h}(Q) < 1.070563363421848 + 2 \hat{h}(G_7)/m^2 < 1.095543 \]
so that
\[ H(Q) < 2.99081. \]
Arguing as in the case of the curve (\ref{ell5}), we must
consider polynomials of the following types, where $a_i \in \Z$:
\begin{eqnarray*}
x^4 + 4 a_1 x^3 + 2 a_2 x^2 + 4 a_3 x + a_4, & |a_1| \leq 2, |a_2| \leq 26, |a_3| \leq 26, |a_4| \leq 80, & \\
x^2+2 a_1 x+a_2, & |a_1| \leq 2, |a_2| \leq 8, & \\
x+a_1, & |a_1| \leq 2. &
\end{eqnarray*}
Each polynomial has to be tested to see if a root can be the $X$-coordinate
of a point in $\mathcal{E}_7(K)$.
Computation shows that in the given range, only the points $\pm G_7$
arise.
It follows that $G_7$ is indeed
a generator of the group of points defined over $K$.


\subsection{The curve $\mathcal{E}_8$ at (\ref{ell8})} \label{E8-computations}
>From the table of Kodaira reductions, we have $\mu_\nu = 0$ except for
\[ \framebox{ $\mu_{\pi}=\frac{1}{4}, \qquad \mu_{\infty_1}=\mu_{\infty_2}=\mu_{\infty_3}=\frac{1}{3}.$ } \]
Further,
\[ \epsilon_\nu^{-1} = \inf_{(X,Y) \in E(K_\nu)} \frac{ \max(|f(X)|_\nu,|g(X)|_\nu)}{\max(1,|X|_\nu)^4} \]
with
\begin{eqnarray*}
f(X) & = & 4 X^3 +(-4-4\ph-2\ph^2-2\ph^3) X^2 + (20+24\ph+6\ph^2+4\ph^3) X, \\
g(X) & = & (X^2 -(5+6\ph+\frac{3}{2}\ph^2+\ph^3))^2.
\end{eqnarray*}
Siksek gives a method for computing the $\epsilon_\nu$. For the
non-Archimedean valuation, we have the following (in Siksek's notation).

First, we observe that 
$g(1+\frac{1}{2}\phi+\frac{1}{4}\phi^3) \equiv 0 \pmod {\pi^{10}}$,
and $g(X) \not \equiv 0 \pmod {\pi^{12}}$ for any $X \in K$. Thus
$\ep_\pi = |\pi|_\pi^{-2 j} = (2^{-\frac{1}{4}})^{-2 j}$, where $j \leq 5$.
This weak inequality is all that we need, resulting in
\[ \framebox{ $\epsilon_{\pi} \leq 2^{\frac{5}{2}}$}. \]

At $\infty_1$,
\[ \epsilon_{\infty_1}^{-1} = \inf_{(X,Y) \in E({\bf R}) } \frac{ \max(|f(X)|,|g(X)|)}{\max(1,|X|)^4}, \]
and the infimum needs to be taken over $X \in {\bf R}$ such that
$f(X) \geq 0$, that is, over
$[0, \infty)$. This infimum occurs at the root $6.700009106939032...$
of $f(X)=g(X)$, and has value $0.52198282519734460776...$,
so that
\[ \framebox{ $\epsilon_{\infty_1} = 1.91577184483403789523...$ }. \]
At $\infty_2$,
\[ \epsilon_{\infty_2}^{-1} = \inf_{(X,Y) \in E({\bf R}) } \frac{ \max(|\bar{f}(X)|,|\bar{g}(X)|)}{\max(1,|X|)^4}, \]
where
\begin{eqnarray*}
\bar{f}(X) & = & 4 X^3 + (-4+4\ph-2\ph^2+2\ph^3) X^2 + (20-24\ph+6\ph^2-4\ph^3) X, \\
\bar{g}(X) & = & (X^2 -(5-6\ph+\frac{3}{2}\ph^2-\ph^3))^2,
\end{eqnarray*}
with infimum taken over $X \in {\bf R}$ such that
$f(X) \geq 0$, that is, over
$[0, \infty)$.
This infimum occurs at the root $0.007079403590926...$ of
$\bar{f}(X)-\bar{g}(X)=0$, and has value $0.00075595704579275884...$
so that
\[ \framebox{ $\epsilon_{\infty_2} = 1322.82648275513803837226...$ }. \]
At $\infty_3$,
\[ \epsilon_{\infty_3}^{-1} = \inf_{(X,Y) \in E({\bf C})} \frac{\max(|F(X)|, |G(X)|)}{\max(1,|X|)^4}, \]
where
\[ F(X)^2 = 16 X^2 (X^4 + (2+\ph^2) X^3 +(8-\ph^2) X^2 + (8+2\ph^2) X +(10+2\ph^2)), \]
\[ G(X) = X^4 + (2+3\ph^2) X^2 + (10+2\ph^2) . \]
The infimum occurs at the root $-0.129793717617543598396...-1.8431948223777410560...i$ of
$|F(z)|=|G(z)|$, with value $0.666554705029609086504527189...$
(on the unit circle, the minimum taken exceeds $8$).
Thus
\[ \framebox{ $\epsilon_{\infty_3} = 1.50025195599748847329235227... $ } \]
Putting the above together results in
\begin{eqnarray*}
h(P)-2 \hat{h}(P) & \leq & \frac{1}{4} (\frac{1}{4} \cdot 4 \cdot \log 2^{\frac{5}{2}} + \frac{1}{3} \cdot 1 \cdot \log(1.91577184483403789523)\\
 & & + \frac{1}{3} \cdot 1 \cdot \log(1322.82648275513803837226) \\
 & & + \frac{1}{3} \cdot 2 \cdot \log(1.5002519559974884732923)),
\end{eqnarray*}
that is,
\[ h(P)-2 \hat{h}(P) \leq 1.153959714852488... \]
Suppose now the point $G_8$ at (\ref{gen8}) is not a generator.
We easily check that $G_8$ is not divisible by $2$ in $E(K)$,
and so $G_8=m Q$ for $m \geq 3$ and $Q \in E(K)$, with
$x(Q) \in \mathcal{O}_K$.
Then
\[ h(Q) \leq 1.153959714852488 + 2 \hat{h}(Q) < 1.153959714852488 + 2 \hat{h}(G_8)/m^2 < 1.167799 \]
so that
\[ H(Q) <  3.21491. \]
Arguing as in the case of the curve (\ref{ell5}), we must
consider polynomials of the following types, where $a_i \in \Z$:
\begin{eqnarray*}
x^4 + 4 a_1 x^3 + 2 a_2 x^2 + 4 a_3 x + a_4, & |a_1| \leq 3, |a_2| \leq 31, |a_3| \leq 33, |a_4| \leq 106, & \\
x^2+2 a_1 x+a_2, & |a_1| \leq 3, |a_2| \leq 10, & \\
x+a_1, & |a_1| \leq 3. &
\end{eqnarray*}
Each polynomial has to be tested to see if a root can be the $X$-coordinate
of a point in $\mathcal{E}_8(K)$.
Computation shows that in the given range, only the points $\pm G_8$,
$\pm G_8+(0,0)$, $\pm 2 G_8+(0,0)$ arise.
It follows that $G_8$ is indeed
a generator of the group of points defined over $K$.


\subsection{The curve $\mathcal{E}_9$ at (\ref{ell9})} \label{E9-computations}
>From the table of Kodaira reductions, we have $\mu_\nu = 0$ except for
\[ \framebox{ $\mu_{\pi}=\frac{1}{4}, \qquad \mu_{\infty_1}=\mu_{\infty_2}=\mu_{\infty_3}=\frac{1}{3}.$ } \]
Further,
\[ \epsilon_\nu^{-1} = \inf_{(X,Y) \in E(K_\nu)} \frac{ \max(|f(X)|_\nu,|g(X)|_\nu)}{\max(1,|X|_\nu)^4} \]
with
\[ f(X) = 4 X^3 +(-8\ph-2\ph^3) X^2 + (4+2\ph^2) X, \quad g(X) = (X^2 -(1+\ph^2/2))^2. \]
Siksek gives a method for computing the $\epsilon_\nu$. For the
non-Archimedean valuation, we have the following (in Siksek's notation).
First, we observe that 
$g(i1+\frac{1}{2}\phi-\frac{1}{4}\phi^3) \equiv 0 \pmod {\pi^{10}}$,
and $g(X) \not \equiv 0 \pmod {\pi^{12}}$ for any $X \in K$. Thus
$\ep_\pi = |\pi|_\pi^{-2 j} = (2^{-\frac{1}{4}})^{-2 j}$, where $j \leq 5$.
This weak inequality is all that we need, resulting in
\[ \framebox{ $\epsilon_{\pi} \leq 2^{\frac{5}{2}}$}. \]

At $\infty_1$,
\[ \epsilon_{\infty_1}^{-1} = \inf_{(X,Y) \in E({\bf R}) } \frac{ \max(|f(X)|,|g(X)|)}{\max(1,|X|)^4}, \]
and the infimum needs to be taken over $X \in {\bf R}$ such that
$f(X) \geq 0$, that is, over
$[0, \infty)$. This infimum occurs at the root $1.432001362205440...$
of $g(X)$, and has value $0.43345064994236763769...$,
so that
\[ \framebox{ $\epsilon_{\infty_1} = 2.30706771378232276809...$ }. \]
At $\infty_2$,
\[ \epsilon_{\infty_2}^{-1} = \inf_{(X,Y) \in E({\bf R}) } \frac{ \max(|\bar{f}(X)|,|\bar{g}(X)|)}{\max(1,|X|)^4}, \]
where
\[ \bar{f}(X) = 4 X^3 + (8\ph+2\ph^3) X^2 + (4+2\ph^2) X, \quad \bar{g}(X) = (X^2 -(1+\ph^2/2))^2, \]
with infimum taken over $X \in {\bf R}$ such that
$f(X) \geq 0$, that is, over
$[0, \infty)$.
This infimum occurs at the root $6.061612256558471...$ of
$\bar{f}(X)-\bar{g}(X)=0$, and has value $0.92450305111791316372...$
so that
\[ \framebox{ $\epsilon_{\infty_2} = 1.08166219547982626230...$ }. \]
At $\infty_3$,
\[ \epsilon_{\infty_3}^{-1} = \inf_{(X,Y) \in E({\bf C})} \frac{\max(|F(X)|, |G(X)|)}{\max(1,|X|)^4}, \]
where
\[ F(X)^2 = 16 X^2 (X^4 -2 X^2 + 2)), \]
\[ G(X) = X^4 +(2+\ph^2) X^2 + 2. \]
The infimum occurs at the root $-4.7565846366129458377743885...i$ of
$|f(z)|=|g(z)|$, with value $0.8788942356277939591822979...$
(on the unit circle, the minimum taken exceeds $4$).
Thus
\[ \framebox{ $\epsilon_{\infty_3} = 1.1377933310550162158769381... $ } \]
Putting the above together results in
\begin{eqnarray*}
h(P)-2 \hat{h}(P) & \leq & \frac{1}{4} (\frac{1}{4} \cdot 4 \cdot \log 2^{\frac{5}{2}} + \frac{1}{3} \cdot 1 \cdot \log(2.30706771378232276809)\\
 & & + \frac{1}{3} \cdot 1 \cdot \log(1.08166219547982626230) \\
 & & + \frac{1}{3} \cdot 2 \cdot \log(1.137793331055016215876938)),
\end{eqnarray*}
that is,
\[ h(P)-2 \hat{h}(P) \leq 0.530938461365339... \]
Suppose now the point $G_9$ at (\ref{gen9}) is not a generator.
We easily check that $G_9$ is not divisible by $2$ in $E(K)$,
and so $G_9=m Q$ for $m \geq 3$ and $Q \in E(K)$, with
$x(Q) \in \mathcal{O}_K$.
Then
\[ h(Q) \leq 0.530938461365339 + 2 \hat{h}(Q) < 0.530938461365339 + 2 \hat{h}(G_9)/m^2 < 0.558878 \]
so that
\[ H(Q) < 1.74871. \]
Arguing as in the case of the curve (\ref{ell5}), we must
consider polynomials of the following types, where $a_i \in \Z$:
\begin{eqnarray*}
x^4 + 4 a_1 x^3 + 2 a_2 x^2 + 4 a_3 x + a_4, & |a_1| \leq 1, |a_2| \leq 9, |a_3| \leq 5, |a_4| \leq 9, & \\
x^2+2 a_1 x+a_2, & |a_1| \leq 1, |a_2| \leq 3, & \\
x+a_1, & |a_1| \leq 1. &
\end{eqnarray*}
Each polynomial has to be tested to see if a root can be the $X$-coordinate
of a point in $\mathcal{E}_9(K)$.
Computation shows that in the given range, only the points $\pm G_9$,
$\pm G_9+(0,0)$, $\pm 2 G_9+(0,0)$ arise.
It follows that $G_9$ is indeed
a generator of the group of points defined over $K$.

\subsection{The curve $\mathcal{E}_{10}$ at (\ref{ell10})} \label{E10-computations}

>From the table of Kodaira reductions, we have $\mu_\nu = 0$ except for
\[ \framebox{ $\mu_{\pi}=\frac{1}{4}, \qquad \mu_{\infty_1}=\mu_{\infty_2}=\mu_{\infty_3}=\frac{1}{3}.$ } \]
Further,
\[ \epsilon_\nu^{-1} = \inf_{(X,Y) \in E(K_\nu)} \frac{ \max(|f(X)|_\nu,|g(X)|_\nu)}{\max(1,|X|_\nu)^4} \]
with
\[ f(X)=  4 X^3 +(-4-2\phi^2) X^2 +(4-2\phi^2) X, \qquad g(X)=(X^2-(1-\frac{1}{2}\phi^2))^2. \]
Siksek gives a method for computing the $\epsilon_\nu$. For the
non-Archimedean valuation, we have the following (in Siksek's notation).
At $\pi$, with $\nu(2)=0$, then
$ \frac{1}{2}\phi+\frac{1}{2}\phi^2-\frac{1}{4}\phi^3 \in U_5 \cap V_5$, 
and $U_6=V_6=\{\}$.
Thus
\[ \framebox{ $\epsilon_{\pi}=|\pi|_{\pi}^{-10}=(2^{-\frac{1}{4}})^{-10}=2^{\frac{5}{2}}$}. \]
At $\infty_1$,
\[ \epsilon_{\infty_1}^{-1} = \inf_{(X,Y) \in E({\bf R}) } \frac{ \max(|f(X)|,|g(X)|)}{\max(1,|X|)^4}, \]
and the infimum needs to be taken over $X \in {\bf R}$ such that
$f(X) \geq 0$, that is, over
$[0, \infty)$. This infimum occurs at the turning point $0.635599759292601...$
of $f(X)=0$, and has value $0.23110328892932097092...$, 
so that
\[ \framebox{ $\epsilon_{\infty_1} = 4.32706953082711459453...$ } \]
At $\infty_2$, since $f$ and $g$ are invariant under $\phi \rightarrow -\phi$, 
we have $\epsilon_{\infty_2} = \epsilon_{\infty_1}$. \\
At $\infty_3$, 
\[ \epsilon_{\infty_3}^{-1} = \inf_{(X,Y) \in E({\bf C})} \frac{\max(|F(X)|, |G(X)|)}{\max(1,|X|)^4}, \]
where
\[ F(X)^2 =  16 X^2 (X^4 + (2+\ph^2) X^3 +(8+\ph^2) X^2 +(8+2\ph^2) X +(10 + 2\ph^2)), \]
\[ G(X) =  X^4 -(6+\ph^2) X^2 +(10+2\ph^2). \]
The infimum occurs at the root $-4.7444841736122500543851896...$ of $|F(z)|=|G(z)|$, 
with value $0.7196560907489582019019193413...$
(on the unit circle, the minimum taken exceeds $10$).
Thus
\[ \framebox{ $\epsilon_{\infty_3} = 1.389552611108012960657952... $ } \]
Putting the above together results in
\begin{eqnarray*}
h(P)-2 \hat{h}(P) & \leq & \frac{1}{4} (\frac{1}{4} \cdot 4 \cdot \log 2^{\frac{5}{2}} + \frac{1}{3} \cdot 1 \cdot \log(4.32706953082711459453)\\
 & & + \frac{1}{3} \cdot 1 \cdot \log(4.32706953082711459453) \\
 & & + \frac{1}{3} \cdot 2 \cdot \log(1.3895526111080129606)),
\end{eqnarray*}
that is,
\[ h(P)-2 \hat{h}(P) \leq 0.732195715015999... \]
Suppose now that the points $P_1$ and $P_2$ at (\ref{gen101}) and
(\ref{gen102}) do not generate the full group of points over $K$.
We first check that that $P_1$ is not divisible in $E(K)$. It is
easy to check that $P_1$ is not divisible by $2$. Suppose $P_1 = m Q$
for $m \geq 3$, for $Q \in E(K)$, with $x(Q) \in \mathcal{O}_K$. Then
\[ h(Q) \leq 2 \hat{h}(Q) + 0.732195715015999 \leq 2/9 \hat{h}(P_1)+ 0.732195715015999 < 0.7532202040..., \]
so that
\[ H(Q) < 2.12383. \] 
If $|\Q(x(Q)):\Q|=4$, then from (\ref{coeffbounds4}), $x(Q)$ is a root
of a polynomial of type
\[ x^4 + 4 a_1 x^3 + 2 a_2 x^2 + 4 a_3 x + a_4, \qquad |a_1| \leq 2, |a_2| \leq 13, |a_3| \leq 9, |a_4| \leq 20; \]
If $|\Q(x(Q)):\Q|=2$, then from (\ref{coeffbounds2}), $X(Q)$ is a root 
of a polynomial of type
\[ x^2 + 2 a_1 x + a_2, \qquad |a_1| \leq 2, |a_2| \leq 4; \]
and if $|\Q(x(Q)):\Q|=1$, then from (\ref{coeffbounds1}), $x(Q)$ is a root
of a polynomial of type
\[ x+ a_1, \qquad |a_1| \leq 2. \]
Search finds that the only points $Q$ satisfying these inequalities
are given by $\pm Q = P_1$, $P_2$, $P_1+(0,0)$, $P_2+(0,0)$,
$P_1 \pm P_2$, $P_1 \pm P_2 +(0,0)$, $2 P_1+(0,0$, and $2 P_2 +(0,0)$.
Since $P_1$ and $P_2$ are of infinite order and independent, it follows
that $P_1$ is not divisible.\\
Further, it
is straightforward to check that the index of the subgroup in $E(K)$
generated by $P_1$ and $P_2$ is {\it odd}.
We take $P_1=G_1$ as one of the generators of $E(K)$, and denote by $G_2$
a second generator. Put $P_2 = a G_1 + m G_2$, for $a, m \in \Z$, 
and where without loss of generality 
\begin{equation}
\label{bds}
m \geq 3, \qquad |a| < m/2. 
\end{equation}
It follows that 
\begin{equation}
m^2 \hat{h}(G_2)=\hat{h}(-a P_1+P_2)=a^2 \hat{h}(P_1)-a <P_1,P_2> +\hat{h}(P_2)
\end{equation}
so that
\[ \hat{h}(G_2) = a^2/m^2 \hat{h}(P_1)-a/m^2 <P_1,P_2>+\hat{h}(P_2)/m^2, \]
whence using (\ref{bds}),
\[ \hat{h}(G_2) < 1/4 \hat{h}(P_1) + 1/6 |<P_1,P_2>| +\hat{h}(P_2)/9 < 0.035009546550940. \]
Thus 
\[ h(G_2) < 2 \hat{h}(G_2) + 0.732195715015999 < 0.8022148081, \]
with 
\[ H(G_2) < 2.23048. \]
As above, $x(G_2)$ is a root of a polynomial of type\\
$\bullet$ $x^4 + 4 a_1 x^3 + 2 a_2 x^2 + 4 a_3 x + a_4, \qquad |a_1| \leq 2, |a_2| \leq 14, |a_3| \leq 11, |a_4| \leq 24$ \\
$\bullet$ $x^2 +2 a_1 x+ a_2, \qquad |a_1| \leq 2, |a_2| \leq 4$ \\
$\bullet$ $x+a_1, \qquad |a_1| \leq 2$.\\
Search finds no points other than those found above in testing
$P_1$ for divisibility, and it follows that indeed $P_1$ and $P_2$
generate the group of points over $K$.

\subsection{The curve $\mathcal{E}_{11}$ at (\ref{ell11})} \label{E11-computations}

>From the table of Kodaira reductions, we have $\mu_\nu = 0$ except for
\[ \framebox{ $\mu_{\pi}=\frac{1}{4}, \qquad \mu_{\infty_1}=\mu_{\infty_2}=\mu_{\infty_3}=\frac{1}{3}.$ } \]
Further,
\[ \epsilon_\nu^{-1} = \inf_{(X,Y) \in E(K_\nu)} \frac{ \max(|f(X)|_\nu,|g(X)|_\nu)}{\max(1,|X|_\nu)^4} \]
with
\begin{eqnarray*}
f(X) & = & 4 X^3 +(-16-20\ph-4\ph^2-4\ph^3) X^2+(52+56\ph+10\ph^2+12\ph^3) X,\\
g(X) & = & (X^2 -(13+14\ph+\frac{5}{2}\ph^2+3\ph^3))^2.
\end{eqnarray*}
Siksek gives a method for computing the $\epsilon_\nu$. For the
non-Archimedean valuation, we have the following (in Siksek's notation).

First, we observe that 
$g(1-\frac{1}{2}\phi+\frac{1}{4}\phi^3) \equiv 0 \pmod {\pi^{10}}$,
and $g(X) \not \equiv 0 \pmod {\pi^{12}}$ for any $X \in K$. Thus
$\ep_\pi = |\pi|_\pi^{-2 j} = (2^{-\frac{1}{4}})^{-2 j}$, where $j \leq 5$.
This weak inequality is all that we need, resulting in
\[ \framebox{ $\epsilon_{\pi} \leq 2^{\frac{5}{2}}$}. \]

At $\infty_1$,
\[ \epsilon_{\infty_1}^{-1} = \inf_{(X,Y) \in E({\bf R}) } \frac{ \max(|f(X)|,|g(X)|)}{\max(1,|X|)^4}, \]
and the infimum needs to be taken over $X \in {\bf R}$ such that
$f(X) \geq 0$, that is, over
$[0, \infty)$. This infimum occurs at the root $6.585832771319615...$
of $f(X)=g(X)$, and has value $0.09399416728471314096...$,
so that
\[ \framebox{ $\epsilon_{\infty_1} = 10.63895802141582727314...$ }. \]
At $\infty_2$,
\[ \epsilon_{\infty_2}^{-1} = \inf_{(X,Y) \in E({\bf R}) } \frac{ \max(|\bar{f}(X)|,|\bar{g}(X)|)}{\max(1,|X|)^4}, \]
where
\begin{eqnarray*}
\bar{f}(X) & = & 4 X^3 + (-16+20\ph-4\ph^2+4\ph^3) X^2 + (52-56\ph+10\ph^2-12\ph^3) X, \\
\bar{g}(X) & = & (X^2 -(13-14\ph+\frac{5}{2}\ph^2-3\ph^3))^2,
\end{eqnarray*}
with infimum taken over $X \in {\bf R}$ such that
$f(X) \geq 0$, that is, over
$[0, \infty)$.
This infimum occurs at the root $0.014878523374045...$ of
$\bar{f}(X)-\bar{g}(X)=0$, and has value $0.00439272523855358661...$
so that
\[ \framebox{ $\epsilon_{\infty_2} = 227.64911204172531486370...$ }. \]
At $\infty_3$,
\[ \epsilon_{\infty_3}^{-1} = \inf_{(X,Y) \in E({\bf C})} \frac{\max(|F(X)|, |G(X)|)}{\max(1,|X|)^4}, \]
where
\[ F(X)^2 = 16 X^2 (X^4 +2\ph^2 X^3 +(6-4\ph^2) X^2 +(4-2\ph^2) X + 2), \]
\[ G(X) = X^4 +(-6+5\ph^2) X^2 + 2. \]
The infimum occurs at the root $-3.839123346088306...$ of
$|F(z)|=|G(z)|$, with value $0.88315461284603293103...$
(on the unit circle, the minimum taken $2.5$).
Thus
\[ \framebox{ $\epsilon_{\infty_3} = 1.13230456530983167167...$ } \]
Putting the above together results in
\begin{eqnarray*}
h(P)-2 \hat{h}(P) & \leq & \frac{1}{4} (\frac{1}{4} \cdot 4 \cdot \log 2^{\frac{5}{2}} + \frac{1}{3} \cdot 1 \cdot \log(10.63895802141582727314)\\
 & & + \frac{1}{3} \cdot 1 \cdot \log(227.64911204172531486370) \\
 & & + \frac{1}{3} \cdot 2 \cdot \log(1.13230456530983167167)),
\end{eqnarray*}
that is,
\[ h(P)-2 \hat{h}(P) \leq 1.103286821056004... \]
Suppose now the point $G_{11}$ at (\ref{gen11}) is not a generator.
We easily check that $G_{11}$ is not divisible by $2$ in $E(K)$,
and so $G_{11}=m Q$ for $m \geq 3$ and $Q \in E(K)$, with
$x(Q) \in \mathcal{O}_K$.
Then
\[ h(Q) \leq 1.103286821056004 + 2 \hat{h}(Q) < 1.103286821056004 + 2 \hat{h}(G_{11})/m^2 < 1.153246 \]
so that
\[ H(Q) < 3.16847. \]
Arguing as in the case of the curve (\ref{ell5}), we must
consider polynomials of the following types, where $a_i \in \Z$:
\begin{eqnarray*}
x^4 + 4 a_1 x^3 + 2 a_2 x^2 + 4 a_3 x + a_4, & |a_1| \leq 3, |a_2| \leq 30, |a_3| \leq 31, |a_4| \leq 100, & \\
x^2+2 a_1 x+a_2, & |a_1| \leq 3, |a_2| \leq 10, & \\
x+a_1, & |a_1| \leq 3. &
\end{eqnarray*}
Each polynomial has to be tested to see if a root can be the $X$-coordinate
of a point in $\mathcal{E}_{11}(K)$.
Computation shows that in the given range, only the points $\pm G_{11}$,
$\pm G_{11}+(0,0)$ arise. It follows that $G_{11}$ is indeed
a generator of the group of points defined over $K$.

\subsection{The curve $\mathcal{E}_{12}$ at (\ref{ell12})} \label{E12-computations}

>From the table of Kodaira reductions, we have $\mu_\nu = 0$ except for
\[ \framebox{ $\mu_{\pi}=\frac{1}{4}, \qquad \mu_{\infty_1}=\mu_{\infty_2}=\mu_{\infty_3}=\frac{1}{3}.$ } \]
Further,
\[ \epsilon_\nu^{-1} = \inf_{(X,Y) \in E(K_\nu)} \frac{ \max(|f(X)|_\nu,|g(X)|_\nu)}{\max(1,|X|_\nu)^4} \]
with
\begin{eqnarray*}
f(X) & = & 4 X^3 +(-12-12\ph-2\ph^2-2\ph^3) X^2 + (20+24\ph+6\ph^2+4\ph^3) X,\\
g(X) & = & (X^2 -(5+6\ph+\frac{3}{2}\ph^2+\ph^3))^2.
\end{eqnarray*}
Siksek gives a method for computing the $\epsilon_\nu$. For the
non-Archimedean valuation, we have the following (in Siksek's notation).

First, we observe that 
$g(1+\frac{1}{2}\phi+\frac{1}{4}\phi^3) \equiv 0 \pmod {\pi^{10}}$,
and $g(X) \not \equiv 0 \pmod {\pi^{12}}$ for any $X \in K$. Thus
$\ep_\pi = |\pi|_\pi^{-2 j} = (2^{-\frac{1}{4}})^{-2 j}$, where $j \leq 5$.
This weak inequality is all that we need, resulting in
\[ \framebox{ $\epsilon_{\pi} \leq 2^{\frac{5}{2}}$}. \]

At $\infty_1$,
\[ \epsilon_{\infty_1}^{-1} = \inf_{(X,Y) \in E({\bf R}) } \frac{ \max(|f(X)|,|g(X)|)}{\max(1,|X|)^4}, \]
and the infimum needs to be taken over $X \in {\bf R}$ such that
$f(X) \geq 0$, that is, over
$[0, \infty)$. This infimum occurs at the turning point
$4.250162195138054...$ of $f(X)/X^4$, and has value $0.14604193738214317332...$,
so that
\[ \framebox{ $\epsilon_{\infty_1} =  6.84734822014400303608...$ }. \]
At $\infty_2$,
\[ \epsilon_{\infty_2}^{-1} = \inf_{(X,Y) \in E({\bf R}) } \frac{ \max(|\bar{f}(X)|,|\bar{g}(X)|)}{\max(1,|X|)^4}, \]
where
\begin{eqnarray*}
\bar{f}(X) & = & 4 X^3 + (-12+12\ph-2\ph^2+2\ph^3) X^2 + (20-24\ph+6\ph^2-4\ph^3) X, \\
\bar{g}(X) & = & (X^2 -(5-6\ph+\frac{3}{2}\ph^2-\ph^3))^2,
\end{eqnarray*}
with infimum taken over $X \in {\bf R}$ such that
$f(X) \geq 0$, that is, over
$[0, \infty)$.
This infimum occurs at the root $0.007463441518832...$ of
$\bar{f}(X)-\bar{g}(X)=0$, and has value $0.00075564996126157299...$
so that
\[ \framebox{ $\epsilon_{\infty_2} = 1323.36405910810826194351...$ }. \]
At $\infty_3$,
\[ \epsilon_{\infty_3}^{-1} = \inf_{(X,Y) \in E({\bf C})} \frac{\max(|F(X)|, |G(X)|)}{\max(1,|X|)^4}, \]
where
\[ F(X)^2 =  16 X^2 (X^4 +(-2+\ph^2) X^3 -3\ph^2 X^2 -4 X +(10+2\ph^2)), \]
\[ G(X) = X^4 +(2+3\ph^2) X^2 +(10+2\ph^2). \]
The infimum occurs at the root $-0.309564888209587...- 5.048223442072423...i$ 
of $|F(z)|=|G(z)|$, with value $0.84342888812084072475...$
(on the unit circle, the minimum taken exceeds $10$).
Thus
\[ \framebox{ $\epsilon_{\infty_3} = 1.18563641118340119834... $ } \]
Putting the above together results in
\begin{eqnarray*}
h(P)-2 \hat{h}(P) & \leq & \frac{1}{4} (\frac{1}{4} \cdot 4 \cdot \log 2^{\frac{5}{2}} + \frac{1}{3} \cdot 1 \cdot \log(6.84734822014400303608)\\
 & & + \frac{1}{3} \cdot 1 \cdot \log(1323.36405910810826194351) \\
 & & + \frac{1}{3} \cdot 2 \cdot \log(1.18563641118340119834)),
\end{eqnarray*}
that is,
\[ h(P)-2 \hat{h}(P) \leq 1.220913082178307... \]
Suppose now the point $G_{12}$ at (\ref{gen12}) is not a generator.
We easily check that $G_{12}$ is not divisible by $2$ in $E(K)$,
and so $G_{12}=m Q$ for $m \geq 3$ and $Q \in E(K)$, with
$x(Q) \in \mathcal{O}_K$.
Then
\[ h(Q) \leq 1.220913082178307 + 2 \hat{h}(Q) < 1.220913082178307 + 2 \hat{h}(G_{12})/m^2 <  1.234753 \]
so that
\[ H(Q) < 3.43753. \]
Arguing as in the case of the curve (\ref{ell5}), we must
consider polynomials of the following types, where $a_i \in \Z$:
\begin{eqnarray*}
x^4 + 4 a_1 x^3 + 2 a_2 x^2 + 4 a_3 x + a_4, & |a_1| \leq 3, |a_2| \leq 35, |a_3| \leq 40, |a_4| \leq 139, & \\
x^2+2 a_1 x+a_2, & |a_1| \leq 3, |a_2| \leq 11, & \\
x+a_1, & |a_1| \leq 3. &
\end{eqnarray*}
Each polynomial has to be tested to see if a root can be the $X$-coordinate
of a point in $\mathcal{E}_{12}(K)$.
Computation shows that in the given range, only the points $\pm G_{12}$,
$\pm G_{12}+(0,0)$, $\pm 2 G_{12}+(0,0)$ arise.
It follows that $G_{12}$ is indeed
a generator of the group of points defined over $K$.

}

\end{document}